\documentclass[a4paper,leqno]{article}
\usepackage[latin1]{inputenc}
\usepackage{amsfonts}
\usepackage{enumitem}
\usepackage{color}
\usepackage{textcomp}
\usepackage [english]{babel}
\usepackage{graphicx,dsfont}
\usepackage{stmaryrd}
\usepackage{amssymb,amsmath,amsthm,}

\newtheorem{define}{Definition}[section]
\newtheorem{pro}{Proposition}[section]
\newtheorem{Lem}{Lemma}[section]
\newtheorem{cor}{Corollary}[section]

\newtheorem{remark}{Remark}[section]
\theoremstyle{plain} \newtheorem{thm}{Theorem}[section]
\newcommand{\eqdef}{\stackrel{\mathrm{def}}{=}}   

\DeclareMathOperator{\dist}{dist}
\catcode`\@=11
\makeatletter

\@addtoreset{equation}{section}
\def\R{\mathbb R}
\def\N{\mathbb N}
\begin{document}
{\vspace{0.001in}}

\title
{Doubly nonlinear equation involving $p(x)$-homogeneous operators: local existence, uniqueness and global behaviour.}
\author{R. Arora, J. Giacomoni and G. Warnault \footnote{LMAP, UMR E2S-UPPA CNRS 5142 B\^atiment IPRA, Avenue de l'Universit\'e F-64013 Pau, France. email: rakesh.arora@univ-pau.fr, jacques.giacomoni@univ-pau.fr guillaume.warnault@univ-pau.fr}}

\maketitle
\begin{abstract}
In this work, we investigate the qualitative properties as uniqueness, regularity and stabilization of the weak solution to the nonlinear parabolic problem involving general $p(x)$-homogeneous operators:
\begin{equation*}
     \left\{
         \begin{alignedat}{2} 
             {} \frac{q}{2q-1}\partial_t(u^{2q-1}) -\nabla.\, a(x, \nabla u)
             & {}= f(x,u) + h(t,x) u^{q-1}  
             && \quad\mbox{ in } \, (0,T) \times \Omega;
             \\
             u & {}> 0
             && \quad\mbox{ in }\, (0,T) \times \Omega ;
             \\
             u & {}= 0
             && \quad\mbox{ on }\, (0,T) \times \partial\Omega;
             \\
             u(0,.)&{}= u_0
             && \quad\  \mbox{in}\, \ \Omega.   
          \end{alignedat}
     \right.
\end{equation*}
Thanks to the Picone's identity obtained in \cite{RAJGGW}, we prove new results about comparison principles which yield a priori estimates, positivity and uniqueness of weak solutions.\vspace{.2cm} \\
\noindent {\textbf{Key words:} Doubly nonlinear equation, Polytropic filtration equation, Leray-Lions operator with variable exponent, stablization.
\vspace{.1cm}
 }
\end{abstract}
\section{Introduction and main results}
The study of various differential equations and variational problems with variable exponent has significantly influenced mathematics in recent years. Indeed, the mathematical problems associated with nonstandard $p(x)$-growth conditions are fascinating in applications as the nonlinear elasticity theory and non-Newtonian fluids models. In specific, the importance of investigating these kinds of problems lies in modeling of various anisotropic features that occur in electrorheological
flows, image restoration, filtration process in complex media, stratigraphy problems and heterogeneous biological interactions.\\   
In the literature, there are many works that explore the questions of existence (local or global), regularity or behaviour of solutions for parabolic equations with variable exponent, for example \cite{AMS, AZ, AS, AS3, RAJGGW,*10, GTW}. Prior investigations have implemented diverse approaches to study the elliptic and parabolic problems with nonstandard growth. In \cite{RAJGGW, GTW}, the authors have followed the semigroup approach, involving the semi-discretization in time method for $p(x)$-Laplacian and Leray-Lions type operators. In \cite{AS}, the Galerkin method is used alternatively to prove the existence of weak solutions and similarly in \cite{AZ, AS3}, the authors have used perturbation methods. We further specify that global properties of solutions as extinction in finite time, localization, blow-up in finite time are also explored in \cite{AS1, AS2}. The existence of mild solutions of parabolic equation and its stabilization properties are studied in \cite{GTW} for $p(x)$-Laplacian and in \cite{*10} for Leray-Lions type operators.\\ 
The original model of our equation is given by 
\begin{equation}\label{basic}
\partial_t u - \nabla.\, (|\nabla (u^m)|^{p-2} \nabla (u^m))=0 \ \ \text{ in }\ (0,T) \times \Omega.
\end{equation}
For $p=2$ and $m>1$, \eqref{basic} is well-known as the porous media equation. More generally, for $p >1$ and $m>0$, \eqref{basic} is known as the Polytropic Filtration Equations (P.F.E.) (see \cite{Wu}). The physical background of P.F.E. can be explained by considering the flow of compressible non-newtonian fluid in the homogeneous isotropic rigid medium which satisfies:
\begin{equation*}
     \left\{
         \begin{alignedat}{2}
             {} \epsilon \partial_t u 
             & {}=  -\nabla (u\overrightarrow{V})\;
             && \quad \mbox{Mass balance}\,
                          \\
             \mathcal{P} & {}= \mathcal{P}_0 u^m
             && \quad\mbox{State equation}\,
          \end{alignedat}
     \right.
\end{equation*} where
$u$ is the particle density of the fluid, $\overrightarrow{V}$ is the momentum velocity, $\mathcal{P}$ is the pressure, $m$ is the polytropic constant and $\mathcal{P}_0$ is the reference pressure and $\epsilon$ is the porosity of the medium. Due to the influence of molecular and ion effects in non-newtonian fluids, the linear Darcy's law is no longer valid. Instead, we have the nonlinear version of Darcy's law:
\begin{equation*}
\mu \overrightarrow{V} = -\lambda |\nabla \mathcal{P}|^{p-1} \nabla \mathcal{P}
\end{equation*}
where $\mu$ is the viscosity of the fluid and $\lambda$ is the permeability of the medium. By combining the two last equations, we obtain an equivalent form of \eqref{basic}.
Depending upon the value of $m$ and $p$, \eqref{basic} is called as Slow Diffusion Equation (S.D.E.) if $p > 1+ \frac{1}{m}$ and Fast Diffusion Equation (F.D.E.) if $p < 1 +\frac{1}{m}$ (for more details see Chapter 2, \cite{Wu}). A main difference between the two cases is the existence of solutions with compact support for the S.D.E whereas the occurrence of dead core type solutions can not occur for the F.D.E. due to the infinite speed of perturbations propagation.\\
In the framework of Doubly Nonlinear Equations (D.N.E. for short) {\it i.e.} $p\neq 2$ and $m\neq 1$, \eqref{basic} is referred in the literature (for instance see \cite{I}) as:
$$\begin{array}{c|c|c}
& p\in(1,2)& p>2\\
\hline
m\in(0,1) & \mbox{Doubly singular} & \mbox{Degenerate-singular}\\
\hline
m>1 & \mbox{Singular-degenerate} & \mbox{Doubly degenerate}\\
\hline
\end{array}$$
The D.N.E. have significant interests because they possess a wide spectrum of applications for instance in fluid dynamics, soil science, combustion theory, reaction chemistry (see \cite{Aris, A, B, B1, bear, childs, GilPet, L, richards, SW} and reference therein) and for D.N.E. involving $p$-Laplacian operator, we refer to \cite{diaz2,kalash, MT, zhan}. The non-homogeneous variant of the model \eqref{basic} together with multivalued source/sink terms can also be interpreted as the limiting case (when $m \to 1$) of the climate Energy Balance Models (see \cite{bedi,berm,diaz}). Recently, the study of D.N.E. involving variable exponent growth are getting into substantial attention: the authors in \cite{RAJGGW} have studied a class of  \eqref{basic} involving the $p(x)$-laplacian using time-discretization method. The authors in \cite{BDMS, BDMS1} and \cite{AS3, AS4} have studied existence of solutions to D.N.E. using a nonlinear version of minimizing movement method and Galerkin method respectively. \\
In the present paper, we study the existence, uniqueness and qualitative properties of the weak solutions of the following D.N.E. driven by a general quasilinear operator of Leray-Lions type:
\begin{equation}\label{Main}
     \left\{
         \begin{alignedat}{2} 
             {} \frac{q}{2q-1}\partial_t(u^{2q-1}) -\nabla.\, a(x, \nabla u)
             & {}= f(x,u) + h(t,x) u^{q-1}  
             && \quad\mbox{ in } \, Q_T;
             \\
             u & {}> 0
             && \quad\mbox{ in }\, Q_T ;
             \\
             u & {}= 0
             && \quad\mbox{ on }\, \Gamma;
             \\
             u(0,.)&{}= u_0
             && \quad\  \mbox{in}\, \ \Omega,   
          \end{alignedat}
     \right.
     \tag{DNE}
\end{equation}
where $T>0$,  $q> 1$, $Q_T \eqdef (0,T) \times \Omega$ with $\Omega \subset \mathbb{R}^N, N \geq 1$ a smooth bounded domain, $\Gamma \eqdef (0,T) \times \partial\Omega$ and $h$ belongs to $L^\infty(Q_T)$.\\
The main difference of this work with the previous studies is the doubly nonlinear feature together combined to the broad class of considered Leray-Lions type operators $a$. More precisely, problem \eqref{Main} involves a class of variational operators $a:\Omega \times \mathbb{R}^N \to \mathbb{R}$ defined as, for any $(x, \xi) \in \Omega \times \mathbb{R}^N$:
$$a(x, \xi)= (a_j(x, \xi))_j\eqdef \left(\frac{1}{p(x)}\partial_{\xi_j} A(x,\xi)\right)_j=\frac{1}{p(x)}\nabla_\xi A(x,\xi)$$
where $A:\Omega \times \mathbb{R}^N \to \mathbb{R^+}$ is continuous, differentiable with respect to $\xi $ and satisfies:
\begin{itemize}
\item[$(A_0)$] $\xi \to A(.,\xi)$ is $p(x)$-homogeneous {\it i.e.} $A(x,t \xi)= t^{p(x)} A(x,\xi),$ for any $t \ \in \mathbb{R}^+$, $\xi \in \mathbb{R}^N $ and a.e. $ x \in \Omega$
\end{itemize}
with $p \in C^1(\overline{\Omega})$ satisfying $$1< p_{-} \eqdef  \min_{x\in \overline{\Omega}} p(x) \leq p(x) \leq p_{+}\eqdef \max_{x\in \overline{\Omega}} p(x) < \infty.$$
This class of operators $a$ also satisfies ellipticity and growth conditions:
\begin{itemize}
        \item[$(A_1)$] For $j\in \llbracket 1,N\rrbracket$,  $a_j(x,0)=0$, $a_j \in C^1(\Omega \times \mathbb{R}^N \backslash \{0\}) \cap C(\Omega \times \mathbb{R}^N)$ and there exist two constants $\gamma,\, \Gamma >0$ such that for all $x \in \Omega$, $\xi \in \mathbb{R}^N \backslash \{0\}$ and $\eta \in 
        \mathbb{R}^N$:
\begin{equation*}
\begin{split}
&\sum_{i,j=1}^N \frac{\partial a_j}{\partial \xi_i}(x, \xi)\ \eta_i \eta_j \geq \gamma |\xi|^{p(x)-2} |\eta|^2;\\
&\sum_{i,j=1}^N \left|\frac{\partial a_j}{\partial \xi_i} (x, \xi) \right| \leq \Gamma |\xi|^{p(x)-2}.
\end{split}
\end{equation*}    
\end{itemize}
\begin{remark}
The assumption $(A_1)$ gives the convexity of $\xi\mapsto A(x,\xi)$ and growth estimates, 
for any $(x,\xi) \in \Omega \times \mathbb{R}^N$: 
\begin{equation}\label{boundA}
\frac{\gamma}{p(x)-1} |\xi|^{p(x)} \leq A(x, \xi) \leq \frac{\Gamma}{p(x)-1} |\xi|^{p(x)}; \qquad |a(x, \xi)| \leq C|\xi|^{p(x)-1}; 
\end{equation} 
and, see \cite{Tolk}, for any $\xi, \eta \in \mathbb{R}^N$ and $x\in \Omega$, there exists a constant $\gamma_0>0$ depending on $\gamma$ and $p$ such that 
\begin{equation}\label{esta}
\langle a(x, \xi)- a(x, \eta), \xi-\eta \rangle \geq
          \gamma_0 \left\{
          \begin{array}{cl}
           |\xi-\eta|^{p(x)} & \mbox{ if }\, p(x)>2;          \\ 
         \displaystyle \frac{|\xi-\eta|^2}{(1+|\xi|+|\eta|)^{2-p(x)}} &\mbox{ if }\, p(x)\leq 2.
          \end{array}
    \right.
\end{equation}
Moreover, the homogeneity assumption implies that $A(x,\xi)=a(x,\xi).\xi$ for any $(x,\xi)\in \Omega\times \R^N$.
\end{remark}
Next, we impose the condition below to insure qualitative properties as regularity and the validity of Hopf Lemma.
\begin{itemize}		  
         \item[$(A_2)$] There exists $C>0$ such that for any $(x,\xi) \in \Omega\times  \mathbb{R}^N \backslash \{0\}$:
         $$\sum_{i,j=1}^N \left|\frac{\partial a_i }{\partial x_j} (x, \xi) \right| \leq 
         C|\xi|^{p(x)-1}(1+|\ln(|\xi|)|).$$ 
\end{itemize}

\begin{remark} 
More precisely, from the condition $(A_2)$ we derive  the Strong Maximum Principle (see \cite{*8}) and the $C^{1, \alpha}$-regularity of weak solutions (see Remark 5.3 in \cite{*2} and Remark 3.1 in \cite{*10}).
\end{remark}
\noindent Concerning the conditions on the functions $f$ and $h$, we assume:
\begin{enumerate}
 \item[$(f_0)$] $f : \overline\Omega \times \mathbb{R^+} \to \mathbb{R^+}$ is a continuous function such that $f(x,0)\equiv 0$ and $f$ is positive on $\Omega\times \R^+\backslash\{0\}$.
 \item[$(f_1)$] For any $x \in \Omega,$ $s\mapsto \frac{f(x,s)}{s^{q-1}}$ is nonincreasing in $\mathbb{R^+}\backslash \{0\}$.
\end{enumerate}
and
\begin{enumerate}
\item[${\bf(H_h)}$] there exists $\underline{h} \in L^{\infty}(\Omega)\backslash\{0\}$, $\underline{h}\geq 0$ such that  $h(t,x) \geq \underline{h} (x)$ for a.e in $Q_T$.
\end{enumerate}
\noindent The study of \eqref{Main} is naturally concerned with the investigation of the following associated parabolic problem:
 \begin{equation}\label{WeakS1}
     \left\{
         \begin{alignedat}{2} 
             {} v^{q-1} \partial_t (v^q) - \nabla.\, a(x, \nabla v)
             & {}= h(t,x) v^{q-1} + f(x,v)
             && \quad\mbox{ in }\, Q_T;
             \\
             v & {}\geq 0
             && \quad\mbox{ in }\, Q_T;
             \\
             v & {}= 0
             && \quad\mbox{ on }\, \Gamma;
             \\
             v(0,.) & {}= v_0           && \quad\mbox{ in }\, \Omega.
          \end{alignedat}
     \right.\tag{E}\end{equation}
We further prove that a weak solution of \eqref{WeakS1} is also a weak solution of \eqref{Main}.\\
By denoting  $\mathcal{W} \eqdef W_0^{1,p(x)}(\Omega)$ (we refer to \cite{*6,KR} for the definitions and properties of variables exponent Lebesgue and Sobolev spaces) and introducing  weighted spaces with the notation $\delta(x)\eqdef\dist(x, \partial \Omega)$:
$$L^\infty_\delta(\Omega)\eqdef\{w: \Omega \to \mathbb{R}\ | \mbox{ measurable},\ \frac{w}{\delta(.)}\in L^\infty(\Omega)\}$$
endowed with the norm $\|w\|_\delta=\sup_\Omega \left|\frac{w(x)}{\delta(x)}\right|$ and for $r>0$: 
\begin{equation*}
\begin{split}
&\mathcal{M}_\delta^{r}({\Omega}) \eqdef  \{ w:\Omega \to \mathbb{R}^+\ |\ \mbox{ measurable}, \ \exists \ c>0,\ \frac{1}{c} \leq \frac{w^{r}}{\delta (x)} \leq c\},
\end{split}
\end{equation*}
we introduce the notion of weak solution of \eqref{WeakS1} as follows:
\begin{define}\label{defw}
Let $T>0$, a weak solution to \eqref{WeakS1} is any positive function $v \in L^{\infty}(0,T;\mathcal{W}) \cap L^\infty(Q_T)$ such that $\partial_t(v^q) \in L^{2}(Q_T)$  satisfying for any $\phi \in  L^2(Q_T) \cap L^1(0,T; \mathcal{W})$  and for any $t \in (0,T]$
\begin{equation}\label{weaksodef}
\begin{split}
\int_0^t \int_{\Omega} \partial_t(v^q) v^{q-1} \phi \,dxds +& \int_0^t \int_{\Omega} a(x,\nabla v). \nabla \phi \,dxds \\
&= \int_0^t \int_{\Omega} (h(s,x) v^{q-1} +f(x,v)) \phi \,dxds
\end{split}
\end{equation} 
and $v(0,.)=v_0$ {\it a.e.} in $\Omega$.
\end{define}
\begin{remark}\label{conversiont}
In sense of Definition \ref{defw}, a solution of \eqref{WeakS1}  belongs to $L^{\infty}(Q_T)$, hence $\frac{q}{2q-1} \partial_t(v^{2q-1}) = v^{q-1} \partial_t(v^q)\in L^2(Q_T)$ holds in weak sense and we deduce the existence of a solution of \eqref{Main}. 
\end{remark}
\vspace{-.3cm}
\begin{remark} 
Prototype examples of operators $a$ satisfying $(A_0)$-$(A_2)$ are given below: for any $(x, \xi) \in \Omega \times \mathbb{R}^N$ and $p\in C^{1,\beta}(\overline{\Omega})$ by:
$$A(x,\xi)=\sum_{j=1}^J\left(g_j(x)\left(\sum_{i\in P_j}\xi_i^2\right)^{\frac{p(x)}{2}}\right)$$ 
where  $(P_j)_{j\in J}$ is a partition of $\llbracket 1,N\rrbracket$, $g_j \in C^{1}(\overline{\Omega}) \cap  C^{0, \beta}(\overline{\Omega})$ and $g_j(x) \geq c>0$ for any $j\in J$.\\
In particular for $A(x,\xi)=|\xi|^{p(x)}$, \eqref{Main} can be classified as S.D.E. if $2q<p_-$  and F.D.E. if $2q>p_+$. 
\end{remark}
\noindent About the existence and properties of solutions of \eqref{WeakS1},  we obtain
\begin{thm}\label{wea}
Let $T>0$ and $q\in(1,p_-)$. Assume $A$ satisfies $(A_0)$-$(A_2)$, $f$ satisfies $(f_0)$, $(f_1)$ and
\begin{itemize}  
\item[$(f_2)$] The mapping $x \mapsto \delta^{1-q}(x)f(x,\delta(x))$ belongs to $L^2(\Omega_\varepsilon)$ for some $\varepsilon>0$ where  $\Omega_\varepsilon\eqdef\{x\in \Omega\, |\, \delta(x)<\varepsilon\}$.
\end{itemize}
\noindent Then, for any $h \in L^{\infty}(Q_T)$ satisfying ${\bf(H_h)}$ and for any initial data $v_0 \in \mathcal{M}_\delta^{1}(\Omega) \cap \mathcal{W}$, there exists a unique solution in sense of Definition \ref{defw}.\\
More precisely, we have:
\item[(i)] Let $v$, $w$ are two weak solutions of \eqref{WeakS1} with respect to the initial data $v_0,\,w_0 \in \mathcal{M}^1_\delta({\Omega}) \cap \mathcal{W}$ and $h$, $g\in L^\infty(Q_T)$ satisfying ${\bf(H_h)}$. Then, for any $t\in [0,T]$:
\begin{equation}\label{conpr}
\Vert (v^q(t)-w^q(t))^+\Vert_{L^2}\leq \Vert (v_0^q-w_0^q)^+\Vert_{L^2}+\int_0^t\Vert (h(s)-g(s))^+\Vert_{L^2}\,ds.
\end{equation}
\item[(ii)] Assume in addition 
A satisfies, for any $x \in \Omega$ and for any $\xi,\,\eta \in \mathbb{R}^N$:
\begin{enumerate}
\item[$(A_3)$] $A(x, \frac{\xi- \eta}{2}) \leq \zeta(x)(A(x,\xi)+A(x,\eta))^{1-s(x)}\left(A(x, \xi)+ A(x, \eta)-2A(x, \frac{\xi+ \eta}{2}) \right)^{s(x)}$\\	where for any $x\in\Omega$,  $s(x)=\min\{1,p(x)/2\}$ and $\zeta(x)=\left({1-2^{1-p(x)}}\right)^{-s(x)}$ 
if  $p(x)<2$ or	$\zeta(x)=\frac12$  if $p(x)\geq2$. 
\end{enumerate}
Then, $v\in C([0,T];\mathcal W)$.
\end{thm} 
\begin{remark}
The above result can be generalized in case $f(x,s)\eqdef\tilde f(x,s)+\hat f(x,s)$ where $\tilde f$ satisfies $(f_1)$ and $s\to \frac{\hat f(x,s)}{s^{q-1}}$ is Lipschitz with respect to the second variable, uniformly in $x\in\Omega$ with constant $\omega>0$. Then if $f$ satisfies additionally $(f_0)$, $(f_2)$ and under same conditions for $A$ and $q$, Theorem \ref{wea} still holds, \eqref{conpr} being replaced by
\begin{equation*}
\Vert (v^q(t)-w^q(t))^+\Vert_{L^2}\leq e^{\omega t}\Vert (v_0^q-w_0^q)^+\Vert_{L^2}+\int_0^t e^{\omega(t-s)}\Vert (h(s)-g(s))^+\Vert_{L^2}\,ds.
\end{equation*}
Similar results have been obtained in \cite{diaz2} in the case of the $p-$laplacian operator.
\end{remark}
\begin{remark}
Prototype example of functions $f$ satisfying $(f_0)$-$(f_2)$ is given by for any $(x,s)\in \Omega\times \R^+$, $f(x,s)=g(x)\delta^\gamma (x)s^\beta$ where $g\in L^\infty(\Omega)$ is a nonnegative function, $\beta\in[0,q-1)$ and $\beta+\gamma > q-\frac32$.
\end{remark}
\begin{remark} The condition $(A_3)$ reformulates the local form of Morawetz-type inequality producing convergence properties.
\end{remark}
\noindent In Theorem \ref{wea}, the uniqueness of the solution in sense of Definition \ref{defw} is obtained by the following theorem relaxing the assumptions on $v_0$ and $h$. More precisely, we show:
\begin{thm}\label{uniqueness}
Let $v,\, w$ are two solutions of \eqref{WeakS1} in sense of Definition \ref{defw} with respect to the initial data $v_0,\,w_0 \in L^{2q}(\Omega)$, $v_0, w_0 \geq 0$ and $h,\,\tilde h\in L^{2}(Q_T)$. Then, for any $t\in[0,T]$:
\begin{equation}\label{maib}
\|v^q(t)- w^q(t)\|_{L^2(\Omega)} \leq \|v_0^q- w_0^q\|_{L^2(\Omega)} + \int_0^t \|h(s)- \tilde{h}(s)\|_{L^2(\Omega)} \,ds. 
\end{equation}
\end{thm}
\noindent Using a similar approach based on nonlinear accretive operators theory as in \cite{BBG,*10,GTW}, we introduce $\mathcal{T}_q: L^2(\Omega) \to L^2(\Omega) $ be the operator with the parameter $q$ defined by 
$$\mathcal{T}_qu= -u^{(1-q)/q}\left(\nabla.\,a(x,\nabla(u^{1/q}))+f(x,u^{1/q})\right)$$
and the associated domain
$$\mathcal{D}(\mathcal{T}_q)= \{w:\Omega \to \mathbb{R}^+\, |\, \mbox{ measurable},\ w^{1/q} \in \mathcal{W}\cap L^{2q}(\Omega),\ \mathcal{T}_qw \in L^2(\Omega)\}.$$
Based on the accretive property of $\mathcal{T}_q$ in $L^2(\Omega)$ (see Theorem \ref{L^2} and Corollary \ref{hidden}) and additional regularity on initial data, we obtain the following stabilization result for the weak solutions to \eqref{WeakS1}:
\begin{thm}\label{stat2}
Under the assumptions of Theorem \ref{wea}, let $v$ be the weak solution of \eqref{WeakS1} with the initial data $v_0 \in \mathcal{M}_\delta^{1}({\Omega}) \cap \mathcal{W}$. Assume that $h\in L^\infty([0,+\infty)\times \Omega)$ satisfying ${\bf(H_h)}$ on $[0,+\infty)\times \Omega$ and 
 there exists $h_\infty\in L^\infty(\Omega)$ such that
\begin{equation}\label{behh}
t^{1+ \eta}\|h(t,.) -h_{\infty}\|_{L^2} =O(1)\  \mbox{ at infinity for some $\eta>0$.}
\end{equation}
Then, for any $r\in[1,\infty)$
$$\|v^q(t,.)- v^q_{stat}\|_{L^{r}} \to 0 \ \text{as} \ t \to \infty$$
where $v_{stat}$ is the unique solution of associated stationary problem with the potential $h_{\infty} \in L^{\infty}(\Omega)$. 
\end{thm}
\begin{remark}
The stabilization in $L^{\infty}$-norm appeals new estimates linked to the $T$-accretivity of the operator $\mathcal{T}_q$ in $L^{\infty}$ and in $ L^1$ (see Remark 1.6  and Theorem 2.1 in \cite{BP} and Theorem 1.18 in \cite{Ha}).
\end{remark}
\begin{remark}
In Theorem \ref{stat2}, we noticed that $v_0 \in \mathcal{M}_\delta^1(\Omega) \cap \mathcal{W}$ implies $v_0^q \in \overline{\mathcal{D}(\mathcal{T}_q)}^{L^2}$ (see Proposition 2.11 in \cite{bb}).
\end{remark}
The current work extends significantly results contained in \cite{RAJGGW} which only apply to the $p(x)$-Laplace operator. We also prove in the present paper new stabilization results for (DNE). To this aim, we borrow the Picone's identity obtained in \cite{RAJGGW} that we recall in the next result:
\begin{thm}[Picone's identity]\label{picone}
Let $B: \Omega \times \mathbb{R}^N \to \mathbb{R}$ be a continuous and differentiable function satisfying $(A_0)$ such that $\xi\mapsto B(x,\xi)$ is strictly convex for any $x\in \Omega$. Let $u,\,v  \in L^\infty(\Omega)$ belonging to $\dot{V}_+^r \eqdef \{w:\Omega \to (0,+\infty) \ | \ w^{\frac1r} \in \mathcal{W}\}$ for  some $r \geq 1$. Then, for any $x\in \Omega$
\begin{equation*}
\frac{1}{p(x)} \nabla_\xi B(x,\nabla (u^{1/r})).\nabla\bigg({\frac{v}{u^{\frac{r-1}{r}}}}\bigg)  \leq B^{\frac{r}{p(x)}}(x,\nabla (v^{1/r})) \  B^{\frac{(p(x)-r)}{p(x)}}(x,\nabla (u^{1/r}))
\end{equation*}
where the inequality is strict if $r >1$ or $\frac{u}{v} \not\equiv \mathrm{const.}>0.$ 
 \end{thm}
\noindent We point out that the general form of operators requires to exploit sharply the Picone's identity. In this regard, the integrability of the quotient in this one forces conditions of regularity or behaviours in the choices of test functions.
\\
Another important part of our work is to study the convergence of the weak solution to a steady state. To this goal, our approach is to use the semigroup theory. Due to the general class of operators, additional  technical computations are needed and performed with the help of  the above Picone's identity.
In our knowledge, the study of solutions of D.N.E involving the class of $p(x)$-homogeneous operators are not discussed so far in literature. So, in this context all results brought in this work are completely new. With both autonomous and non-autonomous terms and the large class of considered operators, \eqref{Main} covers a large spectrum of physical situations. In our study, we also provide new strong maximum principle and weak comparison principle in frame of the large class of operators $a$.\\ 
Turning to the layout of the paper: in section \ref{sec2}, we study a problem related to the parabolic problem \eqref{WeakS1} establishing existence and uniqueness results (Theorem \ref{exis}-\ref{L^2}, Corollary \ref{L^infty}-\ref{hidden}). In section \ref{sec3} we then prove Theorem \ref{wea}. Precisely,  we prove Theorems \ref{existence} concerning the existence of a weak solution in sense of Definition \ref{defw} by semi-discretization in time method. Then the subsection \ref{unik} yields the proof of Theorem \ref{uniqueness} and Corollary \ref{extrauniq} giving the uniqueness using Picone's identity and Theorem \ref{regularity} establishes the regularity of solutions and then completes the proof of Theorem \ref{wea}. Finally in section \ref{sec4}, we establish Theorem \ref{stat2} via a classical argument of semigroup theory. In Appendix \ref{apA}, we state in the framework of general class of operators a weak comparison principle, strong maximum principle and regularity results.   
\section{Elliptic problem related to D.N.E.}\label{sec2}
In this section, we study a class of elliptic problem related to D.N.E. in order to prove Theorem \ref{wea}. First we start with a direct application of Theorem \ref{picone} which provides a comparison principle, uniform estimates and uniqueness.
\begin{Lem}\label{implem}
Let $A: \Omega \times \mathbb{R}^N \to \mathbb{R}$ be a continuous and differentiable function satisfying $(A_0)$ with $a(x, \xi)=\frac{1}{p(x)}\nabla_\xi A(x,\xi)$ such that $\xi\rightarrow A(x,\xi)$ is strictly convex for any $x\in \Omega$. Then, for $r \in [1,p_-)$,  for any $w_1,w_2 \in \mathcal{W} \cap L^\infty(\Omega)$ two positive functions and for any $x\in \Omega$
\begin{equation}\label{implem1}
a(x, \nabla w_1). \nabla\left(\frac{w_1^r-w_2^r}{w_1^{r-1}}\right) + a(x,\nabla w_2). \nabla\left(\frac{w_2^r-w_1^r}{w_2^{r-1}}\right) \geq 0.
\end{equation}
If the equality occurs in \eqref{implem1}, then  $w_1 \equiv w_2$ in $\Omega$.
\end{Lem}
 \begin{proof} Let $w_1,w_2 \in \mathcal{W} \cap L^\infty(\Omega)$ such that $w_1,w_2>0$ in $\Omega$. Then Theorem \ref{picone} yields
 \begin{equation*}
 A^{r/p(x)}(x, \nabla w_1) A^{(p(x)-r)/p(x)}(x, \nabla w_2) \geq a(x,\nabla w_2). \nabla \left(\frac{w_1^r}{w^{r-1}_2}\right).
 \end{equation*}
Then, by using Young inequality and the equality $A(x,\xi)= a(x,\xi) . \xi$, we obtain
 $$a(x, \nabla w_2). \nabla \left( w_2- \frac{w_1^r}{w_2^{r-1}}\right) \geq \frac{r}{p(x)}(A(x, \nabla w_2)-A(x, \nabla w_1)).$$
 Reversing the role of $w_1$ and $w_2$:
 $$a(x, \nabla w_1). \nabla \left( w_1- \frac{w_2^r}{w_1^{r-1}}\right) \geq \frac{r}{p(x)}(A(x, \nabla w_1)-A(x, \nabla w_2))$$
and adding the above inequalities we obtain \eqref{implem1} and the rest of the proof follows from Theorem $5.2$ in \cite{*10}.  
\end{proof}
\subsection{$L^\infty$-potential}\label{sec21}
\noindent In this subsection, we study the following associated elliptic problem:
\begin{equation}\label{E_1}
     \left\{
         \begin{alignedat}{2} 
             {} v^{2q-1} -\lambda \nabla.\, a(x, \nabla v)  
             & {}= h_0 v^{q-1} + \lambda f(x,v)
             && \quad\mbox{ in }\,  \Omega \,;
             \\
             v & {}\geq 0
             && \quad\mbox{ in }\,  \Omega \,;
             \\
             v & {}= 0
             && \quad\mbox{ on }\,  \partial\Omega \, ,
          \end{alignedat}
     \right.
\end{equation}
where $h_0\in L^\infty(\Omega)$ and $\lambda$ is a positive parameter. The notion of weak solution of \eqref{E_1} is defined as follows:
\begin{define}\label{ws}
A weak solution of \eqref{E_1} is any nonnegative and nontrivial function $v \in \mathbf{X}\eqdef  \mathcal{W}\cap L^{2q}(\Omega)$ such that for any  $\phi\in {\mathbf{X}}$
\begin{equation}\label{64}
    \int_{\Omega} v^{2q-1} \phi \,dx + \lambda \int_{\Omega} a(x,\nabla v). \nabla \phi \,dx = \int_{\Omega} h_0 v^{q-1} \phi \,dx + \lambda \int_{\Omega} f(x,v) \phi \,dx.
\end{equation}
\end{define}

\noindent The first theorem gives the existence and the uniqueness of the weak solution of \eqref{E_1}.
\begin{thm}\label{exis}
Assume that $A$ satisfies $(A_0)$-$(A_2)$ and  $f$ satisfies $(f_0)$ and $(f_1)$.
Then, for any  $q\in(1, p_-)$, $\lambda>0$ and $h_0 \in L^{\infty} (\Omega)\backslash\{0\}$, $h_0\geq 0$, there exists a weak solution $v \in C^1(\overline{\Omega})\cap\mathcal{M}_\delta^{1}({\Omega})$ to \eqref{E_1}. \\ Moreover, let $v_1, v_2 $ are two weak solutions to \eqref{E_1} with $h_1,\,h_2\in L^\infty(\Omega)\backslash\{0\}$, $h_1,\,h_2\geq 0$ respectively, we have with the  notation $t^+ \eqdef \max\{0, t\}$:
\begin{equation}\label{Acr}
\|(v_1^q-v_2^q)^+\|_{L^2} \leq \|(h_1-h_2)^+\|_{L^{2}}.
\end{equation}
\end{thm}
\begin{proof}
Define the energy functional $\mathcal{J}$ on $\mathbf{X}$:
\begin{equation}\label{defJ}
\begin{split}
\mathcal{J}(v) =& \frac{1}{2q}\int_{\Omega}v^{2q}\,dx + \lambda \int_{\Omega} \frac{A(x, \nabla v)}{p(x)}\,dx - \frac{1}{q} \int_{\Omega} h_0 (v^+)^q \,dx \\
&- \lambda \int_{\Omega} F(x,v)\,dx
\end{split}
\end{equation}
where $\displaystyle F(x,t)=\int_0^{t^+} f(x,s) ds$.
\\
Note from $(f_0)$-$(f_1)$ that there exists $C>0$ large enough such that for any $(x,s) \in \Omega \times \mathbb{R}^+$
\begin{equation}\label{oo}
 0 \leq f(x,s) \leq C(1+  s^{q-1}).
\end{equation}
By \eqref{boundA} and \eqref{oo},  $\mathcal{J}$ is well defined, continuous on $\mathbf{X}$ 
and we have
\begin{align*}
\mathcal{J}(v)&\geq \|v\|^{q}_{L^{2q}} \bigg(c_1 \|v\|^{q}_{L^{2q}} -c_2\bigg)+  \|v\|_{\mathcal{W}}\left(c_3 \|v\|^{p_- -1}_{\mathcal{W}} - c_4\right)
  \end{align*}
where the constants do not depend on $u$. Thus we deduce that $\mathcal J$ is coercive on $\mathbf{X}$. Therefore we affirm that there exists $v_0 \in\mathbf{X}$ a global minimizer of $\mathcal{J}$.\\
Noting that, with the notation $t^-=t^+-t$,
$$    \mathcal{J}(v_0) \geq \mathcal{J}(v_0^+) + \frac{1}{2q}\int_{\Omega}(v_0^-)^{2q} \,dx +  \lambda \int_{\Omega} \frac{A(x, \nabla v_0^-)}{p(x)} \,dx\geq \mathcal{J}(v^+_0)$$
we deduce $v_0\geq 0$. Let $\phi \in C_c^1(\Omega)$ be a nonnegative and nontrivial function, thus for any $t>0$
$$\mathcal{J}(t\phi) \leq t^q (c_1 t^{q} +c_2 t^{p_--q}- c_3)$$
where the constants are independent of $t$ and $c_3>0$ since $h_0\not\equiv 0$. Hence for $t$ small enough, $\mathcal{J}(t\phi) <0$ and since $\mathcal J(0)=0$, we deduce $v_0\not\equiv 0$. The G\^ateaux differentiability of $\mathcal J$ insures that $v_0$
satisfies \eqref{64}.\\
From Proposition \ref{reg3}, we deduce $v_0\in L^\infty(\Omega)$ and Theorem 1.2 in \cite{*2} provides the $C^{1, \alpha}(\overline{\Omega})$-regularity of $v_0$ for some $\alpha \in (0,1)$.\\
By $(f_0)$ and $(f_1)$, $f$ satisfies $\lim_{s\rightarrow 0^+} f(x,s)s^{1-2q}=\infty$ uniformly in $x\in \Omega$, hence  Lemma \ref{SMPl} implies $v_0 \in \mathcal{M}_\delta^{1}(\Omega)$.\\
Finally, let $v_1, v_2 \in  \mathcal{M}_\delta^{1}({\Omega})$ are two weak solutions of \eqref{E_1} with respect to $h_1$ and $h_2$ respectively. Namely, for any $\phi,\ \Psi \in \mathbf{X}$, we have
$$ \int_{\Omega} v_1^{2q-1} \phi \,dx + \lambda \int_{\Omega} a(x, \nabla v_1). \nabla \phi \,dx = \int_{\Omega} h_1 v_1^{q-1} \phi \,dx + \lambda  \int_{\Omega} f(x,v_1) \phi \,dx$$
and
$$ \int_{\Omega} v_2^{2q-1} \Psi \,dx + \lambda \int_{\Omega} a(x, \nabla v_2). \nabla \Psi \,dx = \int_{\Omega} h_2 v_2^{q-1} \Psi \,dx + \lambda \int_{\Omega} f(x,v_2) \Psi \,dx.$$
Subtracting above expressions by taking $\phi= \left( v_1 - \frac{v_2^q}{v_1^{q-1}}\right)^+$ and $\Psi= \left( v_2 - \frac{v_1^q}{v_2^{q-1}}\right)^-$ then by $(f_1)$ and Lemma \ref{implem}, we obtain
\begin{align*}
\begin{split}
  \int_{\Omega} ((v_1^q- v_2^q)^+)^2 \,dx &\leq \int_{\Omega} (h_1-h_2) (v_1^q - v_2^q)^+ \,dx\\
  & \leq \|(h_1-h_2)^+\|_{L^2(\Omega)} \|(v_1^q-v_2^q)^+\|_{L^2}
  \end{split}
\end{align*}
from which \eqref{Acr} follows.
\end{proof}
\begin{remark}
In the proof of Theorem \ref{exis}, condition $(f_1)$ is not optimal to obtain the existence of a minimizer and to apply Lemma \ref{SMPl}. Indeed define a more general condition on $f$
\begin{itemize}
 \item[$(f'_1)$]  $\limsup_{s  \to +\infty} \frac{f(x,s)}{s^{p_- -1}} < \gamma \Lambda p_{\pm}\  \text{uniformly in}\ x \in \Omega$
\end{itemize}
where $p_{\pm}:= \frac{p_-}{p_+(p_+-1)} $ and $\Lambda^{-1} \eqdef (\sup_{\|u\|_{\mathcal{W}}=1}(\|u\|_{L^{p_-}(\Omega)}))^{p_-}$, condition $(f'_1)$ is a sufficient condition to obtain the existence of a weak solution of \eqref{E_1}. Moreover, to apply Lemma \ref{SMPl} we assume in addition that $f$ satisfies:
\begin{itemize}
 \item[$(f''_1)$]  $\liminf_{s  \to 0^+} \frac{f(x,s)}{s^{2q -1}} > 1\   \text{uniformly in}\ x \in \Omega$.
\end{itemize}
\end{remark}
\begin{remark}\label{runiq}
Inequality \eqref{Acr} implies the uniqueness of the solution in the sense of Definition \ref{ws}. Moreover to obtain \eqref{Acr}, we use more precisely $\phi,\ \psi$ belong to $L^\infty_\delta(\Omega)\cap \mathcal W$. The uniqueness can be also obtained by using Theorem \ref{uniq3}.
\end{remark}
\begin{remark}\label{Racr}
For $q=1$, \eqref{E_1} becomes \begin{equation}\label{q=1}
     \left\{
         \begin{alignedat}{2} 
             {} v+\lambda\mathcal T_1  & {}= h_0
             && \quad\mbox{ in }\,  \Omega \,;
             \\
             v & {}= 0
             && \quad\mbox{ on }\,  \partial\Omega.
          \end{alignedat}
     \right.
\end{equation}
For any $h_0\in L^\infty(\Omega)$ and for any $f\in L^\infty(\overline\Omega\times \R)$ satisfying $(f_1)$ with $q=1$, following the proof of Theorem \ref{exis}, we get the existence of a unique weak solution $v_0\in \mathcal W\cap L^2(\Omega)$ (not necessary nonnegative) in sense of Definition \ref{ws} with $ \phi\in \mathcal W\cap L^2(\Omega)$.\\
Moreover, choosing as test function $\phi=(v_0\pm M)^+$ where $M=\|h_0\|_{L^\infty}+\|f\|_{L^\infty}$, we deduce $v_0\in L^\infty(\Omega)$ and hence for any $\lambda>0$, $R(I+\lambda\mathcal T_1)=L^\infty(\Omega)$.\\
Moreover, let $v_1$ and $v_2$ are two solutions to \eqref{q=1} with $h_1,\, h_2\in L^\infty(\Omega)$ respectively, we get from \eqref{esta} and $(f_1)$: for any $\ell:\R\mapsto\R$ Lipschitz and nondecreasing function such that $\ell(0)=0$:
$$\int_\Omega \mathcal (\mathcal{T}_1v_1-\mathcal T_1v_2)\ell(v_1-v_2)\,dx \geq 0.$$
Thus, by section I.4. in \cite{Ha}, $\mathcal T_1$ is $\,T$-accretive in $L^1(\Omega)$ namely for any $h_1,\,h_2\in L^\infty(\Omega)$ and respectively $v_1,\, v_2$ the solutions to \eqref{q=1},  we have
$$\|(v_1-v_2)^+\|_{L^1} \leq \|(h_1-h_2)^+\|_{L^{1}}.$$
Finally, using Remark 1.6 in \cite{BP}, $\mathcal T_1$ is $T$-accretive in $L^m(\Omega)$, for any $m\in[1,\infty]$ {\it i.e}
$$\|(v_1-v_2)^+\|_{L^m} \leq \|(h_1-h_2)^+\|_{L^{m}},\quad m\in[1,\infty].$$ We point out that $T$-accretivity of $\mathcal T_q$, for $q>1$, in $L^2(\Omega)$ is equivalent to 
$$\int_\Omega \mathcal (\mathcal{T}_qv_1-\mathcal{T}_qv_2)\ell(v_1-v_2)\,dx \geq 0$$ with the fixed choice $\ell(t)=t^+$.
\end{remark}
\vspace{0.2cm}
\noindent In the way of Remark \ref{Racr}, Theorem \ref{exis} implies existence, uniqueness and accretivity results for the perturbed problem induced by the operator $\mathcal{T}_q$:
\begin{cor}\label{L^infty}
Assume $A$ satisfies $(A_0)$-$(A_2)$ and $f$ verifies $(f_0)$ and $(f_1)$. Then, for any $q\in(1, p_-)$, $\lambda>0$ and $h_0 \in L^{\infty} (\Omega)\backslash\{0\}$, $h_0\geq 0$, there exists a unique solution $u \in C^1(\overline{\Omega})$ of
\begin{equation}\label{E_12}
     \left\{
         \begin{alignedat}{2} 
             {} u +\lambda \mathcal{T}_q u
             & {}= h_0 
             && \quad\mbox{ in }\,  \Omega;
             \\
u&>0 && \quad\mbox{ in }\,  \Omega;\\
             u & {}= 0
             && \quad\mbox{ on }\,  \partial\Omega .
          \end{alignedat}
     \right.
\end{equation}
Namely, $u$ belongs to $ \dot V^{q}_+\cap\mathcal M^{1/q}_\delta(\Omega)$ and satisfies: 
\begin{equation}\label{kk}
\int_{\Omega} u \psi \,dx+ \lambda \int_{\Omega} a(x,\nabla (u^{\frac1q})).\nabla(u^{\frac{1-q}{q}}\psi ) - f(x, u^{\frac1q})u^{\frac{1-q}{q}} \psi\,dx = \int_{\Omega} h_0 \psi \,dx
\end{equation}
for any  $\psi$ such that 
\begin{equation}\label{ft}
|\psi|^{1/q}\in L^\infty_\delta(\Omega) \ \mbox{ and }\ \frac{|\nabla \psi|}{\delta^{q-1}(.)} \in L^{p(x)}(\Omega).
\end{equation}
Moreover, if $u_1$ and $u_2$ are two solutions of \eqref{E_12} corresponding to $h_1$ and $h_2$ respectively, then
\begin{equation}\label{accreti}
    \|(u_1-u_2)^+\|_{L^2} \leq \|(u_1-u_2+\lambda (\mathcal{T}_qu_1- \mathcal{T}_qu_2))^+\|_{L^{2}}.
\end{equation}
\end{cor}
\begin{proof}
Define the energy functional $\mathcal{E}$ on $\dot V^{q}_+\cap L^2(\Omega)$ as $\mathcal{E}(u) =\mathcal{J}(u^{1/q})$ where $\mathcal{J}$ is defined in \eqref{defJ}.\\
Let $v_0$ is the weak solution of \eqref{E_1} and the global minimizer of \eqref{defJ}. We set $u_0 = v_0^q$. Then, $u_0$ belongs to $\dot{V}_+^q\cap \mathcal M^{1/q}_\delta(\Omega)$.\\
Let $\psi$ satisfying \eqref{ft}. Then there exists $ t_0> 0$ such that for $t \in (-t_0, t_0)$, $u_0 +t\psi >0$.
Hence we have $\mathcal{E}(u_0 +t \psi)\geq \mathcal{E}(u_0)$ for any $t \in (-t_0, t_0)$. Using Taylor expansion, dividing by $t$ and passing to the limit as $t \to 0$ we deduce that $u_0$ verifies \eqref{kk}.\\
Consider $\tilde u\in \dot{V}_+^q\cap \mathcal M^{1/q}_\delta(\Omega)$ another solution satisfying \eqref{kk}. Thus $\tilde v=\tilde u^{1/q}$ verifies \eqref{64} for $\phi \in L^\infty_\delta(\Omega)\cap \mathcal W$. By Remark \ref{runiq}, we deduce $\tilde v=v_0$ and the uniqueness of the solution of \eqref{E_12}. Finally \eqref{accreti} follows from \eqref{Acr}. 
\end{proof}
\subsection{Extensions for $L^2$-potential}\label{sec22}
\noindent We now generalize existence results of subsection \ref{sec21} for $h_0\in L^2(\Omega)$ by approximation method.
\begin{thm}\label{L^2}
Assume $A$ satisfies $(A_0)$-$(A_2)$ and $f$ verifies $(f_0)$ and $(f_1)$. Then, for any $q\in(1,p_-)$, $\lambda>0$ and $h_0\in L^2(\Omega)\backslash\{0\}$, $h_0\geq0$, there exists a positive weak solution $v\in \mathbf{X}$ of \eqref{E_1} in the sense of Definition \ref{ws}.  Moreover, if $h_0 \in L^r(\Omega)$ for some $r> \max\left\{1, \frac{N}{p_-}\right\}$, $v\in L^\infty(\Omega)$ and $v$ is unique.
\end{thm}
\begin{proof}
Consider $h_n\in C^1_c(\Omega)$, $h_n \geq 0$  which converges to $h$ in $L^2(\Omega)$.  By Theorem \ref{exis}, for any $n\geq 1$, define $v_n\in C^{1, \alpha}(\overline{\Omega})\cap \mathcal M^1_\delta(\Omega)$ as the unique positive weak solution of \eqref{E_1} with $h_0=h_n$.\\
For any $s>1$ and $a, b \geq 0 $, observe that
\begin{equation}\label{alg}
|a-b|^{2s}\leq (a^s-b^s)^2.
\end{equation}
Hence \eqref{Acr} implies, for any $n,\ p\in \N^*$: 
$$\|(v_n-v_p)^+\|_{L^{2q}}\leq \|(v_n^q-v^q_p)^+\|^q_{L^{2}} \leq \|(h_n-h_p)^+\|^q_{L^{2}}.$$
Thus we deduce that $(v_n)$ converges to $v$ in $L^{2q}(\Omega)$ and $(v^q_n)$ converges to $v^q$ in $L^{2}(\Omega)$.\\
Note that the limit $v$ does not depend to the choice of the sequence $(h_n)$ by \eqref{Acr}. So define in particular, for any $n\in \mathbb{N}^*$, $h_n=\min\{h,n\}$. By \eqref{Acr}, we deduce that $(v_n)$ is nondecreasing and for any $n\in \mathbb{N}^*$,
\begin{equation}\label{min}
v(x)\geq v_n(x)\geq v_1(x)\geq c\delta(x)>0 \quad \mbox{{\it a.e.} in  }\Omega,
\end{equation}
for some $c$ independent of $n$.\\
From \eqref{boundA}, \eqref{oo} and using H\"older inequality, equation \eqref{64} with $\phi= v_n$ becomes
\begin{equation*}
\begin{split}
\frac{ \lambda \gamma}{p_+-1} \int_{\Omega} |\nabla v_n|^{p(x)}\,dx & \leq \int_{\Omega} a(x, \nabla v_n). \nabla v_n \,dx \\
&\leq  c \ (\|v_n\|^q_{L^{2q}}(\|h_n\|_{L^2}+1)+\|v_n\|_{L^{2q}})
\end{split}
\end{equation*}
where $c$ does not depend on $n$. Hence we deduce that $(v_n)$ is uniformly bounded in ${\mathcal{W}}$ 
and $v_n$ converges weakly to $v$ in $\mathcal{W}$ (up to a subsequence)
.\\
Now taking $\phi=v_n-v$ in \eqref{64},  we obtain as $n\to \infty$
\begin{equation*}
\left|\int_\Omega f(x,v_n)(v_n-v)\,dx\right|+\left| \int_{\Omega} h_n v_n^{q-1} (v_n-v)\,dx\right| +\left| \int_{\Omega} v_n^{2q-1} (v_n-v)\,dx\right| \to 0
\end{equation*}
which infers $\displaystyle\int_{\Omega} a(x,\nabla v_n) .\nabla (v_n - v)\,dx \to 0$.\\
Since $v_n\rightharpoonup v$ in $\mathcal{W}$, we deduce that:
\begin{equation*}\label{pata}
\int_{\Omega} ( a(x,\nabla v_n) - a(x, \nabla v)).\nabla (v_n - v)\,dx \to 0.
\end{equation*}
Thus we infer that 
 \begin{align}\label{strongconv}
  \int_{\Omega} |\nabla(v_n-v)|^{p(x)} \,dx \to 0 \ \ \text{as} \ \ n \to \infty.  
 \end{align}
Indeed we split $\Omega$ into two parts: $\Omega^{l} =\{x \in \Omega: p(x) \leq 2 \}$ and $\Omega^{u} =\{x \in \Omega: p(x) > 2 \}$.\\
Since $\gamma_0>0$, \eqref{esta} implies \eqref{strongconv} directly on $\Omega^u$. On $\Omega^l$, we get from the H\"older inequality and $(v_n)$ bounded in $\mathcal{W}$:
\begin{equation*}
\begin{split}
    \int_{\Omega^l} |\nabla(v_n-v)|^{p(x)} \,dx &\leq c \left\|\frac{|\nabla(v_n-v)|^{p(x)}}{(|\nabla v|+ |\nabla v_n|)^{r(x)}}\right\|_{L^{\frac{2}{p(x)}}(\Omega^l)} \eqdef c \mathcal{N}\\
&\leq c \left(\int_{\Omega^l} \frac{ |\nabla(v_n-v)|^{2} \,dx}{(|\nabla v|+ |\nabla v_n|)^{2-p(x)}} \,dx \right)^{\hat p}          
\end{split}
\end{equation*}
where  $r(x)= \frac{p(x)(2-p(x))}{2}$, $\hat p=\min\{1,\frac{p_+}{2}\}$ if $ \mathcal{N} \leq 1$ and $\hat p=\frac{p_-}{2}$ otherwise.\\
Hence from \eqref{esta}, we conclude \eqref{strongconv} in $\Omega^l$ and the convergence of $(v_n)$ to $v$ in $\mathcal{W}$.
Then by using dominated convergence Theorem and classical compactness arguments,  we obtain $$a(x, \nabla v_n) \to a(x,\nabla v)\ \ \text{in}\ \ \left(L^{\frac{p(x)}{p(x)-1}(\Omega)}\right)^N.$$
Finally passing to the limit in \eqref{64} satisfied by $v_n$ and applying the dominated convergence Theorem, we obtain $v$ is a weak solution of \eqref{E_1}. 
The regularity arises from Proposition \ref{reg3}.
\end{proof}
\noindent Next result is the extension of Corollary \ref{L^infty} for  $L^2$-potential. 
\begin{cor}\label{hidden}
Assume $A$ satisfies $(A_0)$-$(A_2)$ and $f$ verifies $(f_0)$ and $(f_1)$. Then, for any $q\in(1, p_-)$, $\lambda>0$ and $h_0 \in L^2 (\Omega)\cap L^r(\Omega)\backslash\{0\}$ for some $r> \max\{1, \frac{N}{p_-}\}$, $h_0\geq 0$, there exists a solution  $u$ of \eqref{E_12}. Namely, $u$ belongs to $\dot{V}_+^q\cap L^\infty(\Omega)$ and satisfies  \eqref{kk} for any  $\psi$ verifying \eqref{ft} and there exists $c>0$ such that $u(x)\geq c\delta^{q}(x)$ {\it a.e.} in $\Omega$.
\end{cor}
\begin{proof}
Noting that the existence of a weak solution $v_0 \in L^{\infty}(\Omega)$ of \eqref{E_1} for $h \in L^2(\Omega)$, can be obtained by global minimization method as in Theorem \ref{exis}, we deduce from Theorem \ref{uniq3} that the solution obtained by Theorem \ref{L^2} is a global minimizer.\\
Then we follow the same scheme as the proof of Corollary \ref{L^infty}. We consider the functional energy $\mathcal{E}$ defined on $\dot V^{q}_+\cap L^2(\Omega)$. We set $u_0 = v_0^q$. Then, $u_0$ belongs to $\dot{V}_+^q\cap L^\infty$ and \eqref{min} implies  $u_0(x)\geq  c\delta^{q}(x)$ {\it a.e.} in $\Omega$.\\
Take $\psi$ satisfying \eqref{ft}, then for $t$ small enough, $\mathcal{E}(u_0 +t \psi)\geq \mathcal{E}(u_0)$. From classical arguments, we deduce that $u_0$ verifies \eqref{kk}.
\end{proof}
\section{Parabolic problem related to D.N.E.}\label{sec3}
In this section, we prove Theorems \ref{wea} by dividing the proof into three main steps: existence,  uniqueness and regularity of weak solution. The proof of Theorem \ref{wea} \textit{(i)}  follows from the proof of Theorem $1.5$ in \cite{RAJGGW} using Lemma \ref{implem}, Theorem \ref{L^2} and Corollary \ref{hidden}. Thus we omit the proof.
\subsection{Existence of a weak solution}
In light of Remark \ref{conversiont} and improving Theorem $1.4$ in \cite{RAJGGW} to $p(x)$-homogeneous operator, we consider the problem \eqref{WeakS1}
with {$ v_0 \in \mathcal{M}_\delta^{1}({\Omega}) \cap \mathcal{W}$. 
\begin{thm}\label{existence}
Under the assumptions of Theorem \ref{wea}, there exists a solution $v$ to \eqref{WeakS1} in sense of Definition \ref{defw}. Furthermore $v$ belongs to $C([0,T];L^r(\Omega))$ for any $r\geq 1$ and there exists $C>0$ such that, for any $t \in [0,T]$:
\begin{equation}\label{subsup}
\frac1C\delta(x)\leq v(t,x)\leq C\delta(x) \quad a.\,e. \ \mbox{in } \Omega.
\end{equation}\end{thm}
\begin{proof}
The sketch of the proof is classical and in particular we follow the proof of Theorem $1.4$ in \cite{RAJGGW}. However, for the convenience of the readers, we give the entire proof due to the general form setting of the operator $a$ which requires technical computations. We proceed in several steps:\\
{\bf Step 1:}  Semi-discretization in time of \eqref{WeakS1}\\
Let $n^\star \in \mathbb{N}^*$ and set $\Delta_t =T/n^\star.$ For $n\in \llbracket 0,n^\star\rrbracket$, we define $t_n= n \Delta_t$ and for $ (t,x) \in [t_{n-1}, t_n)\times \Omega$ :
 $$h_{\Delta_t} (t,x)= h^n(x)  \eqdef  \frac{1}{\Delta_t} \int_{t_{n-1}}^{t_n} h(s,x) ds.$$
Thus $\|h_{\Delta_t}\|_{L^{\infty}(Q_T)} \leq \|h\|_{L^\infty(Q_T)}$ and $h_{\Delta_t} \to h $ in $L^2(Q_T).$\\
Applying Theorem \ref{exis} with $\lambda= \Delta_t$, $h_0= \Delta_t h^n + v_{n-1}^q$, we define the implicit Euler scheme, 
\begin{equation}\label{Itersch1}
     \left\{
         \begin{alignedat}{2} 
             {}  \bigg(\frac{v^q_n-v^q_{n-1}}{\Delta_t}\bigg) v_n^{q-1} -\nabla.\, a(x,\nabla v_n) 
             & {}=  h^n v_n^{q-1}+ f(x,v_n)
             && \quad\mbox{ in }\,  \Omega \,;
             \\
             v_n & {}\geq 0
             && \quad\mbox{ in }\,  \Omega \,;
             \\
             v_n & {}= 0
             && \quad\mbox{ on }\,  \partial\Omega \,,
        \end{alignedat}
     \right.
\end{equation}
where, for all $n \in \llbracket 1,n^\star\rrbracket$, $v_n \in C^{1}(\overline{\Omega}) \cap \mathcal{M}_\delta^{1}({\Omega}) $ is the weak solution in sense of Definition \ref{ws} .\\
{\bf Step 2:} Sub- and supersolution \\
In this step, we establish the existence of a subsolution $ \underline{w}$ and a supersolution $\overline{w}$ of a suitable equations such that $v_n\in[\underline{w}, \overline{w}]$ for all $n \in \llbracket 0,n^\star\rrbracket$.\\ 
As in Theorem \ref{exis}, we prove, for any $\mu >0$, there exists a unique weak solution, $\underline{w}_{\mu} \in C^{1}(\overline{\Omega}) \cap \mathcal{M}_\delta^{1}(\Omega)$, to
\begin{equation}\label{subsolu}
     \left\{
         \begin{alignedat}{2} 
             {} - \nabla.\, a(x,\nabla w)  
             & {}=  \mu(\underline{h} w^{q-1}+ f(x,w))
             && \quad\mbox{ in }\,  \Omega \,;
             \\
             w& {}\geq  0
             && \quad\mbox{ in }\,  \Omega \,;             
             \\
             w& {}= 0
             && \quad\mbox{ on }\,  \partial\Omega, \,
        \end{alignedat}
     \right.
\end{equation}
where $\underline{h}$ is defined in ${\bf(H_h)}$.\\
Let $\mu_1 < \mu_2$ and $\underline w_{\mu_1}, \underline w_{\mu_2}$ be weak solutions of \eqref{subsolu}. Then, 
$$\int_{\Omega} a(x, \nabla \underline w_{\mu_1}). \nabla \phi \,dx = \mu_1 \int_{\Omega} (\underline{h}\, \underline w_{\mu_1}^{q-1} + f(x, \underline w_{\mu_1})) \phi \,dx$$
$$\int_{\Omega} a(x, \nabla \underline w_{\mu_2}). \nabla \psi \,dx = \mu_2 \int_{\Omega} (\underline{h}\, \underline w_{\mu_2}^{q-1} + f(x, \underline w_{\mu_2})) \psi \,dx.$$
Summing the above equations with $\phi = \frac{(\underline w_{\mu_1}^q- \underline w_{\mu_2}^q)^+}{\underline w_{\mu_1}^{q-1}}$ and $\psi = \frac{(\underline w_{\mu_2}^q- \underline w_{\mu_1}^q)^-}{\underline w_{\mu_2}^{q-1}}$,
then from \eqref{implem} and $(f_1)$, we deduce $(\underline{w}_{\mu})_\mu$ is nondecreasing. From Theorem 1.2 of \cite{*2} and for some $\mu_0>0$ we obtain, $\|\underline{w}_{\mu} \|_{C^{1}(\overline{\Omega})} \leq C_{\mu_0}$ for any $\mu\leq \mu_0$. Moreover, using Theorem \ref{regi}, we have $\|\underline{w}_\mu\|_{L^{\infty}} \to 0$ as $\mu \to 0$.\\
Therefore $\{\underline{w}_{\mu}: \mu \leq \mu_0\}$ is uniformly bounded and equicontinuous in $C^1(\overline\Omega)$. Applying Arzela-Ascoli Theorem, we obtain, up to a subsequence, $\underline{w}_{\mu} \rightarrow 0$ in $C^1(\overline\Omega)$ as $\mu \rightarrow 0$. Then by Mean Value Theorem, we choose $\mu$ small enough such that  $\underline{w}\eqdef \underline{w}_\mu  \in C^{1}(\overline{\Omega}) \cap \mathcal{M}_\delta^{1}({\Omega})$ satisfies $0< \underline{w} \leq v_0.$  \\
Similarly, there exists $\overline{w}_\kappa \in C^1(\overline\Omega) \cap \mathcal{M}_\delta^{1}({\Omega})$ the weak solution of the following problem:
\begin{equation}\label{supprob}
     \left\{
         \begin{alignedat}{2} 
             {}  -\nabla.\ a(x,\nabla w)  
             & {}= \|h\|_{L^{\infty}(Q_T)} w^{q-1} + f(x,w) + \kappa
             && \quad\mbox{ in }\,  \Omega \,;
             \\
             w & {}\geq 0
             && \quad\mbox{ in }\,  \Omega \,; 
             \\
             w & {}= 0
             && \quad\mbox{ on }\,  \partial\Omega \,.
        \end{alignedat}
     \right.
\end{equation}
By Theorem \ref{regi} and by comparison principle, we have for $\kappa$ large enough that $\overline{w}\eqdef \overline{w}_\kappa\geq w_\kappa\geq v_0$ where $w_\kappa$ is the weak solution of \eqref{lambd}.\\
Rewrite \eqref{Itersch1} as follows
\begin{align*}
v_n^{2q-1} - \Delta_t \nabla.\, a(x,\nabla v_n) =\Delta_t \left(h^nv_n^{q-1} +f(x, v_n) \right)+ v_{n-1}^qv_n^{q-1}.
\end{align*}
Since $\underline{w}\leq v_0 \leq \overline{w}$ and  $\underline{w}$, $\overline{w}$ are respectively a sub- and supersolution of the above equation for $n=1$, Theorem \ref{uniq3} yields  $v_1$ belongs to $[\underline{w},\overline{w}]$ and by induction $v_n\in [\underline{w},\overline{w}]$ for any $n \in \llbracket 1,n^\star\rrbracket.$\\
{\bf Step 3:} A priori estimates \\
Define the functions for $n\in \llbracket 1,n^\star\rrbracket$ and $t\in [t_{n-1}, t_n)$
\begin{align*}
    v_{\Delta_t}(t)= v_n \ \ \text{and} \ \ \tilde{v}_{\Delta_t}(t)= \frac{t-t_{n-1}}{\Delta_t}(v_n^q-v_{n-1}^q)+ v_{n-1}^q
\end{align*}
which satisfy
\begin{align}\label{passto}
    v_{\Delta_t}^{q-1} \partial_t \tilde{v}_{\Delta_t} - \nabla.\, a(x, \nabla v_{\Delta_t})= f(x,v_{\Delta_t}) + h^n v_{\Delta_t}^{q-1}
\end{align} 
and by \textbf{Step 2}, there exists $c>0$ independent of $\Delta_t$ such that for any $(t,x)\in Q_T$
\begin{equation}\label{unifb}
\frac1c \delta(x) \leq v_{\Delta_t},\,\tilde{v}_{\Delta_t}^{1/q}\leq c\delta(x).
\end{equation}
In  \eqref{Itersch1}, summing from $1$ to $n' \in \llbracket 1,n^\star\rrbracket$ and multiplying $\frac{v_n^q-v_{n-1}^q}{v_n^{q-1}}\in \mathbf X$, Young's inequality implies
\begin{equation}\label{n1}
\begin{split}
    \frac{1}{2} \sum_{n=1}^{n'} \int_{\Omega} \Delta_t \bigg(\frac{v_n^q- v_{n-1}^q}{\Delta_t}\bigg)^2 dx&+ \sum_{n=1}^{n'} \int_{\Omega} a(x, \nabla v_n). \nabla\bigg(\frac{v_n^q-v_{n-1}^q}{v_n^{q-1}}\bigg) \,dx \\
    &\leq \sum_{n=1}^{n'} \Delta_t \|h^n\|^2_{L^2} +\sum_{n=1}^{n'} \Delta_t  \left\| \frac{f(x,v_n)}{v_n^{q-1}}\right\|_{L^2}^2.
\end{split}
\end{equation}
Since $v_n\in[\underline{w},\overline{w}]\subset \mathcal M_\delta^1(\Omega)$, \eqref{oo} and $(f_2)$ insure  that $\left(\frac{f(x,v_n)}{v_n^{q-1}}\right)_n$ is uniformly bounded in $L^2(\Omega)$ in $\Delta_t$. Hence, combining \eqref{boundA}, \eqref{n1}, $A(x, \xi)= a(x, \xi). \xi$ and Lemma \ref{implem}, we deduce, for any $n'\geq 1$:
\begin{equation*}
\int_{\Omega}  \frac{c_1|\nabla v_{n'}|^{p(x)}-c_2|\nabla v_0|^{p(x)}}{p(x)}\,dx\leq 
\sum_{n=1}^{n'} \int_{\Omega} a(x,\nabla v_n). \nabla\bigg(\frac{v_n^q-v_{n-1}^q}{v_n^{q-1}}\bigg)dx\leq c_3
\end{equation*}
where the constants are independent of $n'$. The above inequality implies that
\begin{align}\label{iterfuncbound1}
   (v_{\Delta_t})\ \text{is bounded in } L^{\infty}(0,T;\mathcal{W}) \ \text{uniformly in }\ \Delta_t
\end{align}
and from \eqref{n1}, we deduce
\begin{align}\label{deribound}
    (\partial_t \tilde{v}_{\Delta_t}) \ \text{is bounded in} \ L^2(Q_T)\  \text{uniformly in} \ \Delta_t.
\end{align}
Moreover, for $\tilde t= \dfrac{t-t_{n-1}}{\Delta_t}$, we have
\begin{equation*}
\begin{split}
\nabla(\tilde{v}^{\frac1q}_{\Delta_t}) &= \bigg(\tilde t +(1-\tilde t) \bigg(\dfrac{v_{n-1}}{v_n}\bigg)^q \bigg)^{\frac{1-q}{q}} \left(\tilde t \nabla v_n +(1-\tilde t)\left(\dfrac{v_{n-1}}{v_{n}}\right)^{q-1}\nabla v_{n-1} \right).
\end{split} 
\end{equation*}
Hence we deduce from \eqref{iterfuncbound1} and \textbf{Step 2} that
\begin{align}\label{iterfuncbound2}
(\tilde{v}^{1/q}_{\Delta_t}) \ \text{is bounded in } L^{\infty}(0,T;\mathcal{W}) \ \text{uniformly in }\ \Delta_t.
\end{align}
Furthermore using \eqref{alg}, \eqref{deribound} implies
\begin{equation}\label{uniq}
 \sup_{[0,T]} \|\tilde{v}_{\Delta_t}^{1/q}- v_{\Delta_t}\|_{L^{2q}(\Omega)}^{2q} \leq \sup_{[0,T]} \|\tilde{v}_{\Delta_t} - v_{\Delta_t}^q\|^2_{L^2(\Omega)} \leq o_{\Delta_t}(1).      
\end{equation}
Gathering \eqref{iterfuncbound1}-\eqref{uniq}, up to a subsequence, $v_{\Delta_t},\, \tilde{v}^{1/q}_{\Delta_t} \overset{\ast}{\rightharpoonup} v$ in $L^{\infty}(0,T;\mathcal{W})$ as $\Delta_t\rightarrow 0$.\\
From \eqref{unifb} and \eqref{deribound} we deduce that $(\tilde{v}_{\Delta_t})_{\Delta_t}$ is equicontinuous in $C([0,T]; L^r(\Omega))$ for any $r\in[1,+\infty)$. Moreover, from \eqref{alg}, we also deduce that $(\tilde{v}^{1/q}_{\Delta_t})_{\Delta_t}$ is uniformly equicontinuous in $C([0,T]; L^r(\Omega))$ for any $r\in[1,+\infty)$. 
Thus, by Arzela Theorem, we get up to a subsequence that for any $r\in[1,+\infty)$
\begin{align}\label{equi2}
    \tilde{v}_{\Delta_t} \to v^q \ \text{in} \ C([0,T];L^r(\Omega)) \mbox{ and }     v_{\Delta_t} \to v \ \text{in} \ L^{\infty}(0,T;L^r(\Omega)),
\end{align}
hence \eqref{unifb} implies \eqref{subsup}. From \eqref{deribound} and \eqref{equi2}, we obtain
\begin{align}\label{deribound2}
    \partial_t \tilde{v}_{\Delta_t}\to \partial_t( v^q) \ \ \text{in} \ L^2(Q_T).
\end{align}
{\bf Step 4:} $v$ satisfies \eqref{weaksodef}\\
From \eqref{equi2} and \eqref{deribound2}, we have as $\Delta_t \to 0^+$
$$\left|\int_{Q_T} v_{\Delta_t}^{q-1} (v_{\Delta_t}-v)\partial_t \tilde{v}_{\Delta_t}\,dxdt\right|+\left|\int_{Q_T} h^n v_{\Delta_t}^{q-1} (v_{\Delta_t}-v)\,dxdt\right| \to 0 $$
and from $(f_0)$, \eqref{unifb} and \eqref{equi2}, we obtain 
 $$\int_{Q_T}  f(x,v_{\Delta_t}) (v_{\Delta_t}-v)\,dxdt\to 0 \ \text{as} \  \Delta_t \to 0^+.$$
Then, multiplying  \eqref{passto} to $(v_{\Delta_t}-v)$ and passing to the limit, we obtain 
 \begin{equation*}
     \int_{Q_T} a(x,\nabla v_{\Delta_t}).\nabla(v_{\Delta_t}-v) \,dxdt \to 0 \ \text{as} \  \Delta_t \to 0^+.
 \end{equation*}
Since $v_{\Delta_t} \overset{\ast}{\rightharpoonup}  v$ in $L^{\infty}(0,T;\mathcal{W})$ and from the above limit, we conclude
 \begin{equation*}\label{conv2}
 \int_{Q_T} (a(x,\nabla v_{\Delta_t})- a(x,\nabla v)).\nabla(v_{\Delta_t}-v) \,dxdt \to 0 \ \text{as} \  \Delta_t \to 0^+.
 \end{equation*}
By \eqref{esta} and classical compactness arguments, we get 
 \begin{align}\label{maincon}
      a(x, \nabla v_{\Delta_{t}}) \to a(x, \nabla v) \ \ \text{in} \ (L^{p(x)/(p(x)-1)}(Q_T))^N.
 \end{align}
Now, we pass to the limit in \eqref{passto}. First we remark that $(v^{q-1}_{\Delta_t})$ converges to $v^{q-1}$ in $L^2(Q_T)$. Indeed \eqref{alg} and \eqref{uniq}-\eqref{equi2} imply  as  $\Delta_t \to 0$:
\begin{equation*}
\begin{split}
         \|v_{\Delta_t}^{q-1} - v^{q-1}\|^{\frac{2q}{q-1}}_{L^2(Q_T)} &\leq C\int_{Q_T} |v_{\Delta_t}^{q-1}-v^{q-1}|^{\frac{2q}{q-1}}\,dxdt\\
				&\leq C\int_{Q_T} |v_{\Delta_t}^{q}-v^{q}|^2\,dxdt\\
         &\leq C\sup_{[0,T]}\left(\|v_{\Delta_t}^{q} - \tilde{v}_{\Delta_t}\|^2_{L^2} + \|\tilde{v}_{\Delta_t} - v^{q}\|^2_{L^2}\right) \to 0.
      \end{split}
 \end{equation*}
Hence plugging \eqref{deribound} and \textbf{Step 1}, we have in $L^2(Q_T)$:
$$ v_{\Delta_t}^{q-1} \partial_t \tilde{v}_{\Delta_t} \rightarrow  v^{q-1}\partial_t(v^q)\quad \mbox{and}\quad  h_{\Delta_t} v_{\Delta_t}^{q-1}\rightarrow h v^{q-1}.$$
 Thus, we deduce, for any $\phi\in L^2(Q_T)$ as $\Delta_t \to 0^+$:
 \begin{equation}\label{disti1}
\left|\int_{Q_T} \left(v_{\Delta_t}^{q-1}\partial_t \tilde{v}_{\Delta_t}-v^{q-1}\partial_t(v^q)\right)\phi\,dxdt\right|+\left|\int_{Q_T}\left(h_{\Delta_t} v_{\Delta_t}^{q-1}- h v^{q-1}\right) \phi \,dxdt\right|\to 0.
\end{equation}
Furthermore from \eqref{oo} and \eqref{unifb}, $(f(x, v_{\Delta_t})\phi)$ is uniformly bounded in $L^2(Q_T)$ in $\Delta_t$ and by \eqref{equi2} we have $f(x,v_{\Delta_t}) \phi\to f(x,v) \phi$ a.e in $Q_T$ (up to a subsequence). Then, by dominated convergence Theorem we obtain
 \begin{equation}\label{disti3}
     \int_{Q_T} f(x,v_{\Delta_t})\phi\,dxdt\to \int_{Q_T} f(x,v) \phi \,dxdt  \ \ \text{as} \ \ \Delta_t \to 0.
 \end{equation}
Finally gathering \eqref{maincon}-\eqref{disti3}, we conclude that $v$ satisfies \eqref{weaksodef} by passing to the limit in \eqref{passto} for any $\phi \in L^2(Q_T) \cap L^1(0,T;\mathcal{W})$. 
\end{proof}
\subsection{Uniqueness}\label{unik}
\textbf{Proof of Theorem \ref{uniqueness}.} Let $\epsilon \in (0,1)$, we take 
\begin{equation}\label{tf}
\phi = \frac{(v+\epsilon)^q-(w+\epsilon)^q}{(v+\epsilon)^{q-1}}\ \mbox{ and }\ \Psi = \frac{(w+\epsilon)^q- (v+\epsilon)^q }{(w+\epsilon)^{q-1}}
\end{equation}
both belonging to $L^2(Q_T) \cap L^1(0,T; \mathcal{W})$, in 
\begin{equation*}
\begin{split}
\int_0^t \int_{\Omega} \partial_t(v^q) v^{q-1} \phi \,dxds &+ \int_0^t \int_{\Omega} a(x,\nabla v). \nabla \phi \,dxds \\
& =  \int_0^t \int_{\Omega} h(s,x) v^{q-1} \phi \,dxds + \int_0^t \int_{\Omega} f(x,v) \phi \,dxds,\\
\int_0^t \int_{\Omega} \partial_t(w^q) w^{q-1} \psi \,dxds &+ \int_0^t \int_{\Omega} a(x,\nabla w). \nabla \psi \,dxds \\
&= \int_0^t \int_{\Omega} \tilde{h}(s,x) w^{q-1} \psi \,dxds + \int_0^t \int_{\Omega} f(x,w) \psi \,dxds
\end{split}
\end{equation*}
and summing the above equalities, we obtain $\mathbf{I}_\epsilon= \mathbf{J}_\epsilon$ where 
\begin{equation*}
      \begin{aligned}
         \mathbf{I}_\epsilon =& \int_0^t \int_{\Omega} \bigg( \frac{\partial_t(v^q) v^{q-1}}{(v+\epsilon)^{q-1}} - \frac{ \partial_t(w^q) w^{q-1}}{(w+ \epsilon)^{q-1}} \bigg) ((v+\epsilon)^q- (w+ \epsilon)^q) \,dxds\\
         &+ \int_0^t \int_{\Omega} a(x,\nabla (v+\epsilon)). \nabla\bigg(\frac{(v+\epsilon)^q - (w+ \epsilon)^q}{(v+ \epsilon)^{q-1}}\bigg)\,dxds\\
          &+ \int_0^t \int_{\Omega} a(x, \nabla (w+\epsilon)). \nabla\bigg(\frac{(w+\epsilon)^q - (v+ \epsilon)^q}{(w+ \epsilon)^{q-1}}\bigg)\,dxds
          \end{aligned}
\end{equation*}
and
 \begin{equation*}
      \begin{aligned}
          \mathbf{J}_\epsilon &=\int_0^t \int_{\Omega}  \bigg( \frac{h v^{q-1}}{(v+\epsilon)^{q-1}} - \frac{\tilde{h} w^{q-1}}{(w+ \epsilon)^{q-1}} \bigg)((v+\epsilon)^q- (w+ \epsilon)^q) \,dxds\\
          &\ \ \ + \int_0^t \int_{\Omega} \bigg( \frac{f(x,v)}{(v+\epsilon)^{q-1}} - \frac{f(x,w)}{(w+ \epsilon)^{q-1}} \bigg)  ((v+\epsilon)^q- (w+ \epsilon)^q) \,dxds. 
      \end{aligned}
\end{equation*}
First we consider $\mathbf{I}_\epsilon$. Since $\frac{w}{w+\epsilon},\, \frac{v}{v+\epsilon} \leq 1$ and $v,\, w\in L^\infty(Q_T)$, we have
$$
    \bigg|\frac{\partial_t(v^q) v^{q-1}}{(v+\epsilon)^{q-1}} - \frac{\partial_t(w^q) w^{q-1}}{(w+ \epsilon)^{q-1}} \bigg| |(v+\epsilon)^q- (w+ \epsilon)^q|\leq C(|\partial_t(v^q)| + |\partial_t(w^q)|)
$$
where $C$ does not depend on $\epsilon$. Moreover, as $\epsilon \to 0$
$$\bigg( \frac{\partial_t(v^q) v^{q-1}}{(v+\epsilon)^{q-1}} - \frac{ \partial_t(w^q) w^{q-1}}{(w+ \epsilon)^{q-1}} \bigg) ((v+\epsilon)^q- (w+ \epsilon)^q) \to \frac{1}{2} \partial_t (v^q- w^q)^2$$
{\it a.e.} in $Q_T$. Then dominated convergence Theorem and Lemma \ref{implem} give
\begin{equation*}\label{unieq1}
\lim_{\epsilon \to 0}\mathbf{I}_\epsilon \geq \frac{1}{2} \int_0^t \int_{\Omega} \partial_t(v^q- w^q)^2\,dxds.
\end{equation*}
In the same way for $\mathbf{J}_\epsilon$, dominated convergence Theorem implies
\begin{equation*}
\begin{split}
    \int_0^t \int_{\Omega} \bigg(\frac{h v^{q-1}}{(v+\epsilon)^{q-1}} &- \frac{\tilde{h} w^{q-1}}{(w+ \epsilon)^{q-1}} \bigg)((v+\epsilon)^q- (w+ \epsilon)^q)\,dxds\\
    & \to \int_0^t \int_{\Omega} (h-\tilde{h})(v^q- w^q)\,dxds.
    \end{split}
\end{equation*}
Moreover Fatou's Lemma gives
\begin{equation*}
\begin{split}
    &\displaystyle\liminf_{\epsilon \to 0}\int_{0}^t \int_{\Omega} \frac{f(x,v)}{(v+\epsilon)^{q-1}} (w+ \epsilon)^q \,dxds \geq  \int_0^t \int_{\Omega} \frac{f(x,v)}{v^{q-1}} w^q \,dxds, 
\\
   &  \displaystyle\liminf_{\epsilon \to 0} \int_0^t \int_{\Omega} \frac{f(x,w)}{(w+\epsilon)^{q-1}} (v+ \epsilon)^q \,dxds \geq  \int_0^t \int_{\Omega} \frac{f(x,w)}{w^{q-1}} v^q \,dxds.
   \end{split}
\end{equation*}
Hence gathering the three last limits and from $(f_1)$, we obtain 
\begin{equation*}
\liminf_{\epsilon \to 0}\mathbf{J}_\epsilon \leq \int_0^t \int_{\Omega} (h-\tilde{h})(v^q- w^q)\,dxds. 
\end{equation*}
Since $\mathbf{I}_\epsilon=\mathbf{J}_\epsilon$,  we conclude using H\"older inequality that for any $t \in [0,T]$
\begin{equation*}
\frac{1}{2} \int_0^t \int_{\Omega} \partial_t (v^q- w^q)^2\,dxds \leq \int_0^t \|h-\tilde{h}\|_{L^2(\Omega)} \|v^q-w^q\|_{L^2(\Omega)} \,ds
\end{equation*}
and by Gr\"onwall Lemma (Lemma A.4 in \cite{bb}) we deduce \eqref{maib}.
\qed\\

\noindent  Hence we conclude the uniqueness of the solution in sense of Definition \ref{defw} in Theorem \ref{wea}:
\begin{cor}\label{extrauniq}
Let $v$ be a solution of \eqref{WeakS1} in sense of Definition \ref{defw} with the initial data $v_0\in L^{2q}(\Omega)$, $v_0\geq 0$ and $h\in L^{2}(Q_T)$. Then, $v$ is unique.
\end{cor}
\noindent From Theorem \ref{existence} and Corollary \ref{extrauniq}, we deduce the existence result for the parabolic problem involving the operator $\mathcal T_q$:
\begin{thm}\label{conversion}
Under the assumptions of Theorem \ref{wea}, for any $u_0$ such that $u_0^{1/q} \in \mathcal{M}_\delta^{1}(\Omega) \cap \mathcal{W}$, there exists a unique weak solution $u\in L^\infty(Q_T)$ of 
\begin{equation}\label{modifi}
     \left\{
         \begin{alignedat}{2} 
             {} \partial_t u + \mathcal{T}_qu
             & {}= h
             && \quad\mbox{ in }\,   Q_T;
             \\
             u & {}> 0
             && \quad\mbox{ in }\,  Q_T;
             \\
             u & {}= 0
             && \quad\mbox{ on }\,       
             \Gamma;
              \\
             u(0,.) & {}= u_0
             && \quad\mbox{ in }\,       
             \Omega,
          \end{alignedat}
     \right.
\end{equation} 
in the sense that:\begin{itemize}
\item $u^{1/q}$ belongs to $L^\infty(0,T;\mathcal W)$, $\partial_t u\in L^2(Q_T)$;
\item there exists $c>0$ such that for any $t\in [0,T]$, $\frac1c\delta^q(x)\leq u(t,x)\leq c\delta^q(x)$ \textit{a.e.} in $\Omega$;
\item $u$ satisfies, for any $t\in [0,T]$:
\end{itemize}
 \begin{equation}\label{www}
\begin{split}
\int_0^t \int_{\Omega} \partial_t u \psi\,dx&ds + \int_0^t \int_{\Omega} a(x,\nabla u^{1/q}).\nabla(u^{\frac{1-q}{q}}\psi)\,dxds\\
& = \int_{0}^t \int_{\Omega} f(x, u^{1/q})u^{\frac{1-q}{q}} \psi \,dxds + \int_0^t \int_{\Omega} h(s,x) \psi \,dxds,
\end{split}
\end{equation}
for any $\psi$ such that
\begin{equation}\label{eto}
|\psi|^{1/q}\in L^\infty(0,T;L^\infty_{\delta}(\Omega)) \mbox{ and } \frac{|\nabla \psi|}{\delta^{q-1}(\cdot)}\in L^1(0,T;L^{p(x)}(\Omega)).
\end{equation}
Moreover, $u$ belongs to $C([0,T];L^r(\Omega))$ for any $r\in[1,+\infty)$. 
\end{thm}
\begin{proof}
Let $v$ be the weak solution of \eqref{WeakS1} in sense of Definition \ref{defw} obtained by Theorem \ref{existence}. Then, setting in \eqref{weaksodef} $u=v^q$  and choosing $\phi=\frac{\psi}{v^{q-1}}$ with $\psi$ satisfying \eqref{eto}, we get the existence of a solution of \eqref{modifi}.\\
Let us consider the uniqueness issue: let $\tilde u$ be another solution of \eqref{modifi}. We set  $\tilde v= \tilde u^{1/q}$ and taking $\psi=v^{q-1}\phi$ with $\phi \in L^\infty(0,T;L^\infty_{\delta}(\Omega))\cap L^1(0,T;\mathcal W)$ in \eqref{www}, we obtain that $\tilde v$ verifies \eqref{weaksodef} with the additional condition $\phi \in L^\infty(0,T;L^\infty_{\delta}(\Omega))$. Since $v$, $\tilde v$ verify \eqref{subsup}, the test functions defined in \eqref{tf} with $v$ and $\tilde v$ belong to $L^\infty(0,T;L^\infty_{\delta}(\Omega))$. Hence \eqref{maib} holds and we conclude the uniqueness.
\end{proof}
\subsection{ Regularity of weak solution}
\begin{thm}\label{regularity}
Under the assumptions of Theorem \ref{wea}, assume in addition $A$ satisfies $(A_3)$. Then, $v$ the weak solution of \eqref{WeakS1} obtained by Theorem \ref{existence} belongs to $C([0,T];\mathcal{W}).$
\end{thm}
\begin{proof}
The proof is similar as the proof of Theorem 1.1, Step 4 in \cite{*10}. However, the nonlinear term in time implies a specific approach in the computations. Hence for the reader's convenience, we include the complete proof.\\
We have $v \in L^{\infty}(0,T; \mathcal{W})\cap C([0,T];L^{p_-}(\Omega))$ and $p \in C^1(\overline\Omega)$, Theorem 8.4.2 in \cite{*6} yields $\mathcal W\subset L^{p_-}(\Omega)$ with compact embedding. So we deduce $t\mapsto v(t)$ is weakly continuous in $\mathcal W$.\\
Moreover, we consider the mapping $\mathcal{K}(v)=\int_{\Omega} \frac{A(x, \nabla v)}{p(x)}dx$ defined in $\mathcal{W}$. 
The convexity of $A$ implies that $\mathcal{K}$ is weakly lower semicontinuous. Thus for any $t_0 \in [0,T]$, we have
\begin{equation}\label{cont1}
\mathcal{K}(v(t_0)) \leq \liminf_{t \to t_0} \mathcal{K}(v(t)).
\end{equation}
In \eqref{Itersch1}, summing from $n'$ to $n''$ and multiplying by $\frac{v_n^q-v_{n-1}^q}{v_n^{q-1}}\in \mathbf X$, we obtain
\begin{align*}
    \sum_{n=n'}^{n''} \int_{\Omega} &\Delta_t \bigg(\frac{v_n^q- v_{n-1}^q}{\Delta_t}\bigg)^2 \,dx + \sum_{n=n'}^{n''} \int_{\Omega} a(x,\nabla v_n).\nabla\bigg(\frac{v_n^q-v_{n-1}^q}{v_n^{q-1}}\bigg) dx \\
    &=  \sum_{n=n'}^{n''} \int_{\Omega} h^n (v_n^q-v_{n-1}^q) \,dx + \sum_{n=n'}^{n''} \int_{\Omega} \frac{f(x,v_n)}{v_n^{q-1}} (v_n^q-v_{n-1}^q) dx.
\end{align*}
As in {\bf Step 4} of the proof of Theorem \ref{existence}, after using Lemma \ref{implem} we pass to the limit as $n\rightarrow \infty$ and we get: for $t \in [t_0, T]$
\begin{equation}\label{cont4}
\begin{split}
\int_{t_0}^t \int_{\Omega} \partial_t(v^q)^2\,dxds + q \mathcal{K}(v(t)) \leq& \int_{t_0}^t \int_{\Omega} h \partial_t(v^q)\,dxds + q\mathcal{K}(v(t_0)) \\
& + \int_{t_0}^t \int_{\Omega} \frac{f(x,v)}{v^{q-1}} \partial_t(v^q)\,dxds.
 \end{split}
\end{equation}
Taking $\limsup$ in \eqref{cont4} as  $t \to t_0^+$  and by \eqref{cont1} we deduce
\begin{equation*}
\lim_{t \to t_0^+} \mathcal{K}(v(t))= \mathcal{K}(v(t_0))
\end{equation*}
and hence we get the right-continuity of $\mathcal{K}$.\\
Now, for $t>t_0$, let $\eta\in(0,t-t_0)$. We multiply \eqref{WeakS1} by $\displaystyle\tau_\eta v = \frac{v^q(.+\eta,.)-v^q}{\eta v^{q-1}} \in L^2(Q_T) \cap L^1(0,T; \mathcal{W})$ and integrate over $(t_0,t) \times \Omega$ and hence by using Theorem \ref{picone} and Young inequality, we obtain: 
\begin{equation}\label{cont3}
\begin{split}
\int_{t_0}^t \int_{\Omega}  v^{q-1}& \partial_t(v^q)\tau_\eta v\,dxds + \frac{q}{\eta} \int_{t_0}^t \mathcal K(v(s+\eta))- \mathcal K(v(s))\,ds\\
& \geq \int_{t_0}^t \int_{\Omega} h v^{q-1} \tau_\eta v \,dxds + \int_{t_0}^t \int_{\Omega} f(x,v) \tau_\eta v\,dxds.
\end{split}
\end{equation} 
Since $v \in L^{\infty}(0,T; \mathcal{W})$ and $\mathcal K$ is right-continuous in $\mathcal W$, by dominated convergence Theorem, we have as $\eta \to 0^+$ 
$$\frac1\eta\int_{t_0}^{t_0+\eta} \mathcal K (v(s))\,ds \to\mathcal K (v(t_0))\quad \mbox{and}\quad \frac1\eta\int_{t}^{t+\eta} \mathcal K (v(s))\,ds \to\mathcal K (v(t)).$$
Then \eqref{cont3} yields,
\begin{equation*}
\begin{split}
\int_{t_0}^t \int_{\Omega} \partial_t(v^q)^2 \,dxds + q\mathcal K (v(t)) \geq& \int_{t_0}^t \int_{\Omega} h \partial_t(v^q)\,dxds+ q\mathcal K (v(t_0))\\
& + \int_{t_0}^t \int_{\Omega} \frac{f(x,v)}{v^{q-1}} \partial_t(v^q)\,dxds.
 \end{split}
\end{equation*}
From \eqref{cont4}, we have the equality for any $t$, $t_0\in[0,T]$ in the above inequality and we deduce the left-continuity of $\mathcal K$.\\
By $(A_3)$, the proof of corollary $A.3$ in \cite{GV} holds by considering $\mathcal K$ as the semimodular. Then, we deduce that $\nabla v(t)$ converges to $\nabla v(t_0)$ in $L^{p(x)}(\Omega)^N$ as $t\to t_0$ and hence $v \in C([0,T]; \mathcal{W}).$
\end{proof}

\section{Stabilization}\label{sec4}
\subsection{Stationary problem related to \eqref{WeakS1}}\label{section0}
In the aim of studying the behaviour of global solution of the problem \eqref{WeakS1} as $t\to \infty$, we consider the following problem
\begin{equation}\label{Stprob}
     \left\{
         \begin{alignedat}{2} 
             {} -\nabla.\, a(x, \nabla v)  
             & {}= b(x)v^{q-1} + f(x,v)
             && \quad\mbox{ in }\,  \Omega;
             \\
             v & {}\geq 0
             && \quad\mbox{ in }\,  \Omega;
             \\
             v & {}= 0
             && \quad\mbox{ on }\,  \partial\Omega,
          \end{alignedat}
     \right.\tag{$S$}
\end{equation}
where $b\in L^\infty(\Omega)$. The notion of weak solution of \eqref{Stprob} is defined as follows:
\begin{define}\label{stdef}
A weak solution to \eqref{Stprob} is any nonnegative function $v \in \mathcal{W}\cap L^{\infty}(\Omega)$, $v\not\equiv 0$ such that for any $\phi \in \mathcal{W}$, $v$ satisfies
\begin{equation}\label{weaksodef11}
\begin{split}
 \int_{\Omega} a(x,\nabla v). \nabla \phi\,dx = \int_{\Omega} b v^{q-1} \phi\,dx + \int_{\Omega} f(x,v) \phi\,dx.
\end{split}
\end{equation} 
\end{define}
\begin{thm}\label{stat1}
Assume that $A$ satisfies $(A_0)$-$(A_2)$ and $(f_0)$ and $(f_1)$ hold. Then, for any $q \in (1, p_-)$, $b\in L^{\infty}(\Omega)\backslash \{0\}$, $b\geq 0$, there exists a unique weak solution $v \in C^{1}(\overline{\Omega})\cap \mathcal{M}_\delta^{1}(\Omega)$ to \eqref{Stprob}.
\end{thm}
\begin{proof}
Consider the energy functional $\mathcal{L}$ defined on $\mathcal{W}$ such that
 \begin{equation*}
\mathcal{\tilde L}(v) =  \int_{\Omega} \frac{A(x, \nabla v)}{p(x)}\,dx - \frac{1}{q} \int_{\Omega}b(v^+)^q\,dx -  \int_{\Omega} F(x,v)\,dx 
\end{equation*}
where $F$ is defined as in \eqref{defJ}. By following the same arguments as in Theorem \ref{exis}, we deduce the existence of nonnegative global minimizer $v_0$ to $\mathcal{L}$ and  the G\^ateaux differentiability of $\mathcal{\tilde L}$ implies $v_0$ satisfies \eqref{weaksodef11}.\\
Combining Proposition \ref{rg1} and Theorem 4.1 in \cite{*3} , we deduce $v_0 \in L^{\infty}(\Omega)$. Then by Theorem 1.2 of \cite{*2}, we obtain, $v_0 \in C^{1, \alpha}(\overline{\Omega})$ for some $\alpha \in (0,1)$. From Lemma \ref{SMPl}, we deduce $v_0 > 0$ and $v_0$ belongs to $\mathcal{M}^1_\delta({\Omega})$.\\
Let $\tilde v_0$ another solution of \eqref{Stprob}. As previously, we deduce that $\tilde v_0\in C^{1,\alpha}(\overline \Omega)\cap \mathcal{M}^1_\delta({\Omega})$. \\
We choose $\displaystyle{\frac{v_0^q-\tilde v_0^q}{v_0^{q-1}}}$ and $\displaystyle\frac{\tilde v_0^q-v_0^q}{\tilde v_0^{q-1}}$ as test functions in \eqref{weaksodef11} satisfied by $v_0$ respectively $\tilde v_0$, then adding the both equations we deduce from Lemma \ref{implem} and $(f_1)$:
\begin{equation*}
\int_{\Omega} a(x, \nabla v_0). \nabla\left(\frac{v_0^q-\tilde v_0^q}{v_0^{q-1}}\right) + a(x,\nabla \tilde{v}_0). \nabla\left(\frac{\tilde v_0^q-v_0^q}{\tilde v_0^{q-1}}\right)\,dx=0.
\end{equation*}
Applying once again Lemma \ref{implem}, we obtain $v_0=\tilde v_0$.  
\end{proof}
\noindent Hence we obtain using the same way of the proof of Corollary \ref{L^infty}:
\begin{cor}
Under the conditions of Theorem \ref{stat1}, there exists a unique solution $u$ of the following problem
\begin{equation}\label{modifi1}
     \left\{
         \begin{alignedat}{2} 
             {}  \mathcal{T}_qu  
             & {}= b
             && \quad\mbox{ in }\,   \Omega;
             \\
             u & {}> 0
             && \quad\mbox{ in }\,  \Omega;
             \\
             u & {}= 0
             && \quad\mbox{ on }\,       
             \partial \Omega.
         \end{alignedat}
     \right.
\end{equation} 
Namely, $u$ belongs to $\dot{V}_+^q\cap \mathcal M^{1/q}_\delta(\Omega)$ and satisfies, for any $\psi$ such that \eqref{ft}:
\begin{equation*}
\int_{\Omega} a(x,\nabla u^{1/q}).\nabla(u^{\frac{1-q}{q}}\psi)\,dx - \int_{\Omega} \frac{f(x, u^{1/q})}{u^{(q-1)/q}} \psi\,dx = \int_{\Omega} b\psi\,dx .
\end{equation*} 
\end{cor}
\subsection{Proof of Theorem \ref{stat2}}

\noindent \textbf{Proof of Theorem \ref{stat2}}.
We consider two cases:\\
\textbf{Case 1:} $h\equiv h_\infty$.\\
We introduce the family $\{S(t); t \geq 0\}$ on $\mathcal{M}_\delta^{1/q}(\Omega)\cap \dot{V}^q_+$ defined as $w(t)=S(t) w_0$ where $w$ is the solution obtained by Theorem \ref{conversion} (and Theorem \ref{existence}) of 
\begin{equation}\label{evol}
     \left\{
         \begin{array}{l l} 
            \partial_t w + \mathcal{T}_qw
             = h_\infty &\mbox{ in }   Q_T;\\
             w > 0 &\mbox{ in }  Q_T;\\
             w = 0 &\mbox{ on }    \Gamma;\\
						w(0,.)=w_0 &\mbox{ in }\Omega.
                     \end{array}
     \right.
\end{equation}
Thus $\{S(t); t \geq 0\}$ defines a semigroup on $\mathcal{M}_\delta^{1/q}(\Omega)\cap\dot{V}^q_+$. Indeed the uniqueness and properties of solution of \eqref{modifi} imply for any $w_0$, 
\begin{equation}\label{semiprop}
S(t+s) w_0 = S(t) S(s) w_0 , \ \   S(0)w_0= w_0
\end{equation} and from \eqref{equi2} the map $t \to S(t)w_0$ is continuous from $[0, \infty)$ to $L^2(\Omega)$.\\
Note that $v=(S(t)w_0)^{1/q}$ is the solution of \eqref{WeakS1} in the sense of Definition \ref{defw} with $h=h_\infty$ and the initial data $w_0^{1/q}$.\\
Let $T>0$ and $v$ be the solution of \eqref{WeakS1} obtained by Theorem \ref{existence} with  $h\equiv h_\infty$ and the initial data $v_0$, hence we get $u(t)=v(t)^q=S(t)u_0$ with $u_0=v_0^q$. \\
Let $\underline w=w_\mu$ be the solution of \eqref{subsolu} and $\overline w=\overline w_\kappa$ be the solution of \eqref{supprob}.  Then, $\underline w,\,\overline w\in\mathcal M_\delta^1(\Omega)$ and for $\mu$ small enough and $\kappa$ large enough, $\underline w$ is a subsolution and $\overline w$ a supersolution of \eqref{Stprob} with $b=h_\infty$ such that $\underline w\leq v_0\leq\overline w$.\\
We define $\underline u(t)=S(t)\underline w^q$ and $\overline u(t)=S(t)\overline w^q$  the solutions to \eqref{evol}. So $\underline u$ and $\overline u$ are obtained by the iterative scheme \eqref{Itersch1} with $v_0=\underline w$ and $v_0=\overline w$. Hence, by construction  the map $t\to \underline u(t)$ is nondecreasing, the map $t\to \overline u(t)$ is nonincreasing and \eqref{conpr} insures for any $t\geq 0$,
\begin{equation}\label{comp}
\underline w^q\leq \underline u(t)\leq u(t)\leq \overline u(t)\leq \overline w^q \ a.\,e. \mbox{ in } \Omega.
\end{equation}
We set $\underline u_\infty=\lim_{t\to \infty} \underline u(t)$ and $\overline u_\infty=\lim_{t\to \infty} \overline u(t)$. Then from \eqref{semiprop}, the continuity in $L^2(\Omega)$ and monotone convergence theorem, we get in $L^2(\Omega)$:
$$\underline u_\infty = \lim_{s \to \infty} S(t+s)(\underline{w}^q) =  S(t)(\lim_{s \to \infty} S(s) (\underline{w}^q))= S(t)\underline u_\infty $$  
and analogously we have $ \overline u_\infty= S(t) \overline u_\infty$. We deduce $\underline u_\infty$ and $\overline u_\infty$ are solutions of \eqref{modifi1} with $b=h_\infty$ and by uniqueness, we have $u_{\text{stat}}\eqdef\underline u_\infty= \overline u_\infty$ where $u_{\text{stat}}$ is the stationary solution of perturbed parabolic problem \eqref{modifi1}.  Therefore from \eqref{comp} and dominated convergence Theorem, we obtain
$$\|u(t)- u_{\text{stat}}\|_{L^2}\to 0 \  \text{as} \ t \to \infty.$$
Finally, using \eqref{comp} and interpolation inequality $\|.\|_r \leq \|.\|_{\infty}^\theta \|.\|_2^{1-	\theta}$, we  conclude the above convergence for any $r\ge 1 $.\\
\textbf{Case 2:} $h\not\equiv h_\infty$.\\
From \eqref{behh}, for any $\varepsilon$ and for some $\eta'\in (0,\eta)$, there exists $t_0>0$ large enough such that for any $t\geq t_0$:
$$ t^{1+\eta'}\|h(t,.)-h_\infty\|_{L^2} \leq \varepsilon.$$
Let $T>0$ and $v$ be the solution of \eqref{WeakS1} obtained by Theorem \ref{existence} with  $h$ and the initial data $v_0=u_0^{1/q}$ and we set $u=v^q$. \\
Since $v$ satisfies \eqref{subsup}, we can define $\tilde u(t)=S(t+t_0)u_0=S(t)u(t_0)$. Then, by \eqref{conpr} and uniqueness, we have for any $t>0$:
\begin{equation*}
\|u(t+t_0,.)-\tilde u(t,.)\|_{L^2}\leq \int_0^{t} \|h(s+t_0,.)-h_\infty\|_{L^2}\,ds\leq \frac{\varepsilon}{t_0^{\eta'}}\leq \varepsilon.
\end{equation*}
By \textbf{Case 1}, we have $\tilde u(t)\to u_{\text{stat}}$ in $L^2(\Omega)$ as $t \to \infty$. Therefore, we obtain
$$\|u(t)- u_{\text{stat}}\|_{L^2}\to 0 \  \text{as} \ t \to \infty$$ 
and by using interpolation inequality we conclude the proof of Theorem \ref{stat2}.
\qed
\appendix 
\section{Additional results}\label{apA}
In this section, we give extensions of technical results for the class of operator $A$ or for some boundary value problems. \\
We begin by extending Theorem 4.3 in \cite{RAJGGW} using Lemma \ref{implem}. Then, we obtain the comparison principle:
\begin{thm}\label{uniq3}
Assume $A$ satisfies $(A_0)$-$(A_2)$ and $f$ satisfies $(f_0)$ and $(f_2)$. 
Let $\underline v,\, \overline{v}\in \mathbf X\cap L^\infty(\Omega) $ be nonnegative functions respectively subsolution and supersolution to \eqref{E_1} for some $h \in L^{r}(\Omega),\ r \geq 2$, $h \geq 0$.  Then $\underline v\leq\overline{v}.$
\end{thm}
\noindent The proof is similar as the proof of Theorem \ref{uniqueness} where the sub- and supersolution do not need to belong to $\mathcal M^1_\delta(\Omega)$. The proof is very similar and we omit it. In the next theorem, we extend Lemma $2.1$ of \cite{*31} and Lemma $3.2$ of \cite{*10} for $p(x)$-homogeneous operators. 
\begin{thm}\label{regi}
Assume $A$ satisfies $(A_0)$-$(A_2)$. Let $\lambda>0$ and  $w_\lambda \in \mathcal{W} \cap C^{1, \alpha}(\overline\Omega)$ be the positive weak solution of  \begin{equation}\label{lambd}
     \left\{
         \begin{alignedat}{2} 
             {} - \nabla.\, a(x,w_\lambda)
             & {}=  \lambda
             && \quad\mbox{ in }\,  \Omega;
             \\
             w_{\lambda} & {}= 0
             && \quad\mbox{ on }\,  \partial\Omega.
        \end{alignedat}
     \right.
\end{equation}
Then, there exists $\lambda^*>0$ such that $w_\lambda$ satisfies
\begin{itemize}
\item for any $\lambda\geq \lambda^*$, $\|w_\lambda\|_{L^{\infty}(\Omega)} \leq C_1 \lambda^{1/(p_--1)}$ and $w_\lambda(x) \geq C_2 \lambda^{\frac{1}{p_+-1+ \varepsilon}} \delta(x)$ for some $\varepsilon \in (0,1)$;
\item for $\lambda < \lambda^*$, $\|w_{\lambda}\|_{L^{\infty}(\Omega)} \leq C_3 \lambda^{1/(p_+-1)}$
\end{itemize}
\noindent where the constants depend upon $p_+, p_-, N,\ \Omega$ and $\alpha$. Moreover if $\lambda_1 < \lambda_2$ then $w_{\lambda_1} \leq w_{\lambda_2}.$ 
\end{thm}
\noindent Now we state a Strong and Hopf maximum principle for variable exponent $p(x)$-homogeneous operators and theirs proof follows from Lemma 3.3 and 3.4 in \cite{RAJGGW}. 

\begin{Lem}\label{SMPl}
Let $\alpha,\, \beta$ be two measurable functions such that $1<\beta_-\leq\beta_+<\alpha_-\leq\alpha_+<\infty$. Let $h,\, l \in L^{\infty}(\Omega)$ be nonnegative functions, $h >0$ and $k:\Omega\times \R^+\to \R^+$ and $A$ satisfies $(A_0)$-$(A_1)$. Consider $u \in C^1(\overline\Omega)$ a nonnegative and nontrivial solution to
\begin{equation*}
     \left\{
         \begin{alignedat}{2} 
             {} -\nabla.\, a(x, \nabla u)  +l(x) u^{\alpha(x)-1} 
             & {}= h(x) u^{\beta(x)-1} + k(x,u)
             && \quad\mbox{ in }\, \Omega \,;
             \\
             u & {}= 0
             && \quad\mbox{ on }\, \partial\Omega\,.
          \end{alignedat}
     \right.
\end{equation*}
If  $\displaystyle\liminf_{t \to 0^+} k(x,t)t^{1-\alpha(x)} > \|l\|_{L^{\infty}}$ uniformly in $x \in \Omega$, then $u$ is positive in $\Omega$.\\
Furthermore, if $\Omega$ satisfies the interior ball condition for any $x\in\partial \Omega$, then $\frac{\partial u}{\partial \vec n} (x)<0$ where $\vec n$ is the outward unit normal vector at $x$.
\end{Lem}
\noindent We state a slight extension of Proposition $A.1$ in \cite{RAJGGW} and Proposition A.2 in \cite{*10}.
\begin{pro}\label{rg1}
Let $q \in [1,p_-)$. Assume $A$ satisfies $(A_0)$-$(A_2)$ and $u\in \mathbf X$ satisfying for any $\Psi \in\mathbf X$:
\begin{equation*}
\int_\Omega a(x,\nabla u) . \nabla \Psi \,dx=\int_\Omega h u^{q-1} \Psi \,dx
\end{equation*}
where $h\in L^2(\Omega)\cap L^{r}(\Omega)$ with $r >\max\{1,\frac{N}{p_-}\}$. Then $u\in L^\infty(\Omega)$.
\end{pro}
\begin{pro}\label{reg3}
Under the assumptions of Proposition \ref{rg1}, consider $u\in \mathbf X$ a nonnegative function satisfying, for any $\Psi \in \mathbf X$, $\Psi\geq 0$:
\begin{equation*}
\int_\Omega u^{2q-1}\Psi \,dx+\int_\Omega a(x,\nabla u) \cdot\nabla \Psi \,dx\leq\int_\Omega(f(x,u)+hu^{q-1}) \Psi \,dx 
\end{equation*}
where $f$ verifies for any $(x,t)\in \Omega\times \R^+$, $|f(x,t)|\leq c_1+c_2|t|^{s(x)-1}$ with $s\in C(\overline\Omega)$ such that for any $x\in\overline \Omega$, $1<s(x)<p^*(x)$ and $h\in L^2(\Omega)\cap L^{r}(\Omega)$ with $r>\max\{1,\frac{N}{p_-}\}$. Then $u\in L^\infty(\Omega)$.
\end{pro}
\noindent The proofs of above results follow the proofs of Theorem 4.1 in \cite{*3} and Proposition A.1 in \cite{RAJGGW}.


\begin{thebibliography}{100}
\bibitem{AMS} E. Acerbi, G. Mingione and G.A. Seregin, {\it Regularity results for parabolic systems related to a class of non-Newtonian fluids}, Ann. Inst. H. Poincare  Anal. Non Lineaire, {\bf 21} (2004),  no. 1, 25-60.
\bibitem{AZ} Y. A. Alkhutov and V. V. Zhikov, {\it Existence theorems for Solutions of Parabolic Equations with Variable Order of Nonlinearity}, Proc. Steklov Inst. Math., {\bf 270} (2010),  no. 1, 15-26.
\bibitem{AS} S. N. Antontsev and S. I. Shmarev, {\it Anisotropic parabolic equations with variable nonlinearity}, Publ. Mat., {\bf 53} (2009), no. 2, 355-399.
\bibitem{AS1} S. N. Antontsev, S. I. Shmarev, {\it A model porous medium equation with variable exponent of nonlinearity: existence, uniqueness and localization properties of solutions}, Nonlinear Anal., {\bf 60} (2005), no. 3, 515-545.
\bibitem{AS2} S. N. Antontsev and S. I. Shmarev, {\it Localization of solutions of anisotropic parabolic equations}, Nonlinear Anal., {\bf 71} (2009), no. 12, 725-737.
\bibitem{AS3} S. N. Antontsev and S. I. Shmarev, {\it Parabolic equations with double variable nonlinearities},  Math. Comput. Simulation, {\bf 81} (2011), no. 10,  2018-2032.
\bibitem{AS4} S. N. Antontsev, S. I. Shmarev, {\it Existence and uniqueness for doubly nonlinear parabolic equations with nonstandard growth conditions}, Differ. Equ. Appl., {\bf 4} (2012), no. 1, 67-94.
\bibitem{Aris} R. Aris, The mathematical theory of diffusion and reaction in permeable catalysts, I, II. Clarendon: Oxford, 1975.
\bibitem{A} D. G. Aronson,  {\it Regularity properties of flows through porous media: The interface}, Arch. Rational Mech. Anal., {\bf 37} (1970), 1-10.
\bibitem{RAJGGW} R. Arora, J. Giacomoni, G. Warnault, {\it A Picone identity for variable exponent operators and applications}, Adv. Nonlinear Anal. {\bf 9} (2020), no. 1, 327-360.
 \bibitem{BBG} M. Badra, K. Bal and J. Giacomoni, {\it A singular parabolic equation: Existence and stabilization},  J. Differential Equations, {\bf 252} (2012), no. 9, 5042-5075.
\bibitem{B} G. I. Barenblatt, {\it  On some unsteady motions of a liquid and gas in a porous medium} (Russian), Akad. Nauk SSSR. Prikl. Mat. Meh., {\bf 16} (1952), 67-78.
\bibitem{B1} G. I. Barenblatt, {\it  On self-similar solutions of the Cauchy problem for a nonlinear parabolic equation of unsteady filtration of a gas in a porous medium} (Russian), Prikl. Mat.
Meh., {\bf 20} (1956), 761-763.
\bibitem{bear} J. Bear, Dynamics of fluids in porous media, New York: Elsevier, 1972.
\bibitem{bedi} S. Bensid and J.I. D\`iaz, {\it On the exact number of monotone solutions of a simplified Budyko climate model and their different stability}, Discrete and Continuous Dynamical Systems, Series {\bf 24} (2019), 1033-1047.
\bibitem{berm} R. Bermejo, J. Carpio, J. I. D\'iaz, L. Tello, {\it Mathematical and numerical analysis of a nonlinear diffusive climate energy balance model}, Math. Comput. Modelling, {\bf 49} (2009), no. 5-6, 1180-1210.
\bibitem{BP} P. Benilan and C. Picard, {\it Quelques aspects non lin\'eaires du principe du maximum} S\'eminaire de Th\'eorie du Potentiel, No. 4, Paris, 1977/1978, Lecture Notes in Math., vol. 713, Springer, Berlin, 1979, pp. 1-37.
\bibitem{BDMS} V. B\"ogelein, F. Duzaar, P. Marcellini and C. Scheven, {\it A variational approach to doubly nonlinear equations},  Atti Accad. Naz. Lincei Rend. Lincei Mat. Appl., {\bf 29} (2018), no. 4, 739-773.
\bibitem{BDMS1} V. B\"ogelein, F. Duzaar, P. Marcellini and C. Scheven, {\it Doubly nonlinear equations of porous medium type}, Arch. Ration. Mech. Anal., {\bf 229} (2018), no. 2, 503-545.
\bibitem{bb} H. Brezis, {\it Operateurs maximaux monotones et semigroupes de contractions dans les espaces de Hilbert}, North Holland, Amsterdam, 1973.
\bibitem{childs} E. C. Childs, An introduction to the physical basis of soil water phenomena, London: Wiley, 1969.
\bibitem{diaz} J.I. D\'iaz, Diffusive Energy Balance Models in Climatology, Stud. Math. Appl., vol. 31, North-Holland, Amsterdam, 2002.
\bibitem{diaz2} J.I. D\'iaz, {\it New applications of monotonicity methods to a class of non-monotone parabolic quasilinear sub-homogeneous problems}, to appear in Pure and Applied Functional Analysis.
\bibitem{*6}  L. Diening, P. Harjulehto, P. H\"ast\"o and M. R{\r u}{\u z}i{\u c}ka, Lebesgue and Sobolev Spaces with Variable Exponents, Springer Berlin Heidelberg, 2011.
\bibitem{*31} X. Fan, {\it On the sub-supersolution method for p(x)-Laplacian equations}, J. Math. Anal. Appl., {\bf 330} (2007), no. 1, 665-682.
\bibitem{*2} X. Fan, {\it Global $C^{1,\alpha}$ regularity for variable exponent elliptic equations in divergence form}, J. Differential Equations, {\bf 235} (2007), no. 2, 397-415
\bibitem{*3} X. Fan and D. Zhao, {\it A class of De Giorgi type and H\"older continuity}, Nonlinear Anal., {\bf 36} (1999), no. 3, Ser. A: Theory Methods, 295-318.
\bibitem{GV} J. Giacomoni and G. Vallet, { \it Some results about an anisotropic $p(x)-$ Laplace-Barenblatt equation}, Adv. Nonlinear Anal., {\bf 1} (2012), no. 3, 277-298.
\bibitem{*10} J. Giacomoni, V. R{\u a}dulescu  and G. Warnault, {\it Quasilinear parabolic problem with variable exponent: Qualitative analysis and stablization}, Commun. Contemp. Math., {\bf 20} (2018), no. 8, 38 pp.
\bibitem{GTW} J. Giacomoni, S. Tiwari and G. Warnault, {\it Quasilinear parabolic problem with p(x)-Laplacian: existence, uniqueness of weak solutions and stabilization}, NoDEA Nonlinear Differential Equations Appl., {\bf 23} (2016), no. 3, Art. 24, 30 pp.
\bibitem{GilPet} B. H. Gilding, L. A. Peletier, {\it The Cauchy problem for an equation in the theory of infiltration}, Arch. Rational Mech. Anal., {\bf 61} (1976), no. 2, 127-140.
\bibitem{Ha} K. S. Ha, {\it Sur les semi-groupes non lin\'eaires dans les espaces $L^\infty(\Omega)$}, J. Math. Soc. Japan, {\bf 31} (1979), no. 4, 593-622.  
\bibitem{I} A. V. Ivanov, {\it Regularity for doubly nonlinear parabolic equations}, J. Math. Sci., {\bf 83} (1997), no. 1, 22-37.
\bibitem{kalash} A. S. Kalashnikov, {\it Some problems of the qualitative theory of nonlinear degenerate second-order parabolic equations}, Uspekhi Mat. Nauk, {\bf 42} (1987), no. 2, 169-222.
\bibitem{L} O. A. Ladyz{\u e}nskaja, {\it New equations for the description of the motions of viscous incompressible fluids, and global solvability for their boundary value problems} (Russian), Trudy Mat. Inst. Steklov., {\bf 102} (1967), 85-104.
\bibitem{KR} V. R{\u a}dulescu and D. Repov\v{s}, Partial Differential Equations with Variable Exponents, Variational Methods and Qualitative Analysis, Monographs and Research Notes in Mathematics, CRC Press, Taylor \& Francis Group, 2015.
\bibitem{richards} L. A. Richards, {\it Capillary conduction of liquids through porous mediums}, Physics, {\bf 1} (1931), no. 5, 318-333.
\bibitem{SW} R. Showalter, N. J. Walkington, {\it Diffusion of fluid in a fissured medium with microstructure}, SIAM J. Math. Anal., {\bf 22} (1991), no. 6, 1702-1722.
\bibitem{Tolk} P. Tolksdorf, {\it Regularity for a more general class of quasilinear elliptic equations}, J. Differential Equations, {\bf 51} (1984), no. 1, 126-150.
\bibitem{MT} M. Tsutsumi, {\it On solutions of some doubly nonlinear degenerate parabolic equations with absorption}, J. Math. Anal. Appl., {\bf 132} (1988), no. 1, 187-212.
\bibitem{Wu} Z. Wu, J. Zhao, J. Yin, H. Li, Nonlinear Diffusion Equations, World Scientific, Singapore, (2001).
\bibitem{zhan} H. Zhan, {\it Infiltration equation with degeneracy on the boundary}, Acta Appl. Math., {\bf 153} (2018), 147-161.
\bibitem{*8} Q. Zhang, {\it A strong maximum principle for differential equations with nonstandard $p(x)$-growth conditions}, J. Math. Anal. Appl., {\bf 312} (2005), no. 1, 24-32.
\end{thebibliography}
\end{document}